\documentclass[a4paper,12pt]{article}

\usepackage[margin=27mm]{geometry}
\usepackage[T1]{fontenc}
\usepackage{amsmath,amssymb,amsthm,mathtools,bm}
\usepackage{booktabs,array,enumitem}
\usepackage{graphicx,float}
\usepackage{subcaption}
\usepackage{xcolor}
\usepackage{hyperref,xurl,enumitem,bm,esint}
\newcommand{\R}{\mathbb R}
\newcommand{\Om}{\Omega}
\newcommand{\Ast}{A^{\mathrm{st}}}

\newcommand{\ip}[2]{\langle #1,#2\rangle}
\newcommand{\norm}[1]{\left\|#1\right\|}
\newcommand{\Tcal}{\mathcal T}

\newcommand{\Leta}{\mathcal L_\eta}

\newcommand{\bb}{\mathbf b}
\newcommand{\bc}{\mathbf c}

\newcommand{\br}{\mathbf r}

\newcommand{\bx}{\mathbf{x}}

\newcommand{\dd}{\,\mathrm d}

\newcommand{\diver}{\nabla \cdot}

\newcommand{\tavg}[1]{\frac1{NT}\sum_{i=1}^{N}\int_0^T #1\dd t}
\newcommand{\tint}[1]{\sum_{i=1}^{N}\int_0^T #1\dd t}
\DeclareMathOperator*{\argmin}{arg\,min}

\theoremstyle{plain}
\newtheorem{theorem}{Theorem}[section]
\newtheorem{lemma}[theorem]{Lemma}
\newtheorem{corollary}[theorem]{Corollary}
\newtheorem{proposition}[theorem]{Proposition}
\theoremstyle{definition}

\theoremstyle{remark}
\newtheorem{remark}[theorem]{Remark}

\numberwithin{equation}{section}

\definecolor{softblue}{RGB}{236,244,251}
\definecolor{softgreen}{RGB}{237,248,241}

\title{\textbf{Model-Based Learning of Urban Crime Dynamics, Part I: identifiability and information loss}\\[0.35em]
}
\author{Baoli Hao$^1$, Annalisa Quaini$^2$, and Ming Zhong$^3$}
\date{}

\begin{document}
\maketitle

\begin{center}
$^1$ Department of Applied Mathematics, Illinois Institute of Technology, 10 West 35th St,
Chicago, IL 60616, USA \\ {\tt bhao2@hawk.illinoistech.edu}

\vskip .2cm
\noindent $^2$ Chair of Computational Methods in Civil and Environmental Engineering, Technical University of Munich,
Arcisstra{\ss}e 21, 80333 Munich, Germany \\ {\tt annalisa.quaini@tum.de}

\vskip .2cm
\noindent $^3$ Department of Mathematics, University of Houston, 3551 Cullen Blvd, Houston, TX 77204, USA \\
{\tt  mzhong3@central.uh.edu}
\end{center}

\begin{abstract}
We consider a partial differential equation (PDE) model of residential
burglary that describes the coupled evolution of the criminal density and
environmental attractiveness.  The model features three spatially varying
parameters (the strength of the neighborhood effects, the static
attractiveness, and a source term for criminal regeneration) that control
the formation of localized crime hotspots.  Since these parameters cannot be
measured directly, we study their recovery from observations of the
attractiveness and criminal density along one or more trajectories.  We propose a weak-form learning method that reduces the
recovery to a linear least-squares problem in a discrete dual norm. We prove that all three
parameters are identifiable under a pointwise rank condition on differences
of the observed attractiveness, and that this condition guarantees the unique
solvability of the learner at the PDE level.  When the static attractiveness is
unknown, we quantify exactly the information loss about the neighborhood effects.  Near a smooth no-flux boundary, the information in the normal
direction decays quadratically with the distance to the boundary, and the smallest singular value of the rank matrix decays linearly.  Experiments with synthetic data 
recover all three parameters with
relative errors of a few percent and illustrate the predicted losses of
information.
\end{abstract}

\noindent \emph{Key words}: 
Parameter identification; 
Weak-form learning;
Crime modeling;
Identifiability;
No-flux boundary conditions

\section{Introduction}
%
%
Predictability and control of urban crime are fundamental for public safety, yet crime remains notoriously hard to predict. Criminological and sociological studies have documented how offenders are driven by environmental cues and respond to the success or failure of their past activities \cite{brantingham1984patterns, wilson1982broken,johnson1997new, wright1994burglars}, and how crime concentrates in a small number of places called hotspots \cite{weisburd2015law}. 
Crime hotspots emerge, migrate, split, and dissipate, with lifetimes and sizes that depend on geographic, economic, and seasonal factors, as observed, e.g., in burglary data from Houston and Chicago \cite{Hao2026,Hao2026b}. Purely data-driven approaches, ranging from deep learning \cite{wang2017deep, wang2019deep} to stochastic and point-process models \cite{cai2022multiscale, mohler2020impact, park2021investigating}, have made substantial progress in forecasting. See, e.g., \cite{blumstein2002crime, dorsogna2015statistical, ersoz2025artificial}. However, they offer limited insight into the mechanisms that generate the observed patterns.

To understand these mechanisms, a probabilistic agent-based model of residential burglary was developed in \cite{short2008statistical}
and the associated PDE system was derived
through a mean-field limit. 
Residential burglary provides a relatively simple setting for such modeling, as mobile offenders target stationary sites. The PDE system in \cite{short2008statistical} 
reproduces the formation of crime hotspots.  This model has since been analyzed and extended in many directions: hotspot formation, dissipation, and displacement \cite{short2010dissipation}, bifurcations and suppression of hotspots \cite{short2010nonlinear}, local existence and uniqueness of solutions \cite{rodriguez2010local}, spike solutions \cite{mei2020existence}, nonlinear diffusion \cite{rodriguez2020relaxation}, L\'evy-flight offender motion \cite{pan2018crime}, incorporation of housing data \cite{woodworth2014nonlocal}, realistic settings \cite{Hao2026}, and law enforcement \cite{jones2010statistical, zipkin2014cops,nurhan2026modelling,Hao2026b}. 

The PDE model of \cite{short2008statistical}, recalled in Sec.~\ref{sec:pd}, depends strongly on three parameters~\cite{Hao2026}. The diffusivity $\eta$, related to near-repeat victimization and the broken-window effect, controls whether hotspots form, their number and size~\cite{short2008statistical,Hao2026}, and the onset of oscillations when police respond with a delay~\cite{Hao2026b}. The static attractiveness $\Ast$ encodes the baseline appeal of residences and anchors persistent hotspots, e.g., along highways~\cite{calamunci2022highway,Hao2026}. The source $q$ determines the homogeneous equilibrium, influences its stability~\cite{short2008statistical}, and affects the average crime level~\cite{Hao2026b}. These parameters are usually prescribed as constants in space, but real cities are heterogeneous, and $\eta$, $\Ast$, and $q$ cannot be measured directly: they must be inferred from data.

Inferring the parameters raises two issues.  The first is identifiability: whether
direct observations of $A$ and $\rho$ along trajectories  determine $\eta$, $\Ast$, and $q$ uniquely.  The second is
computation: how to recover the parameters accurately and efficiently.
In this paper, we address both.  We propose a
weak-form learning method in two variants: one learns $\eta$ and $q$ when
$\Ast$ is known, and the other learns $\eta$, $\Ast$, and $q$
simultaneously.  We will show that the
three parameters are identifiable simultaneously, and that the discrete
learner is uniquely solvable for sufficiently rich test spaces. 

The main contributions of this paper are the following: i) \emph{Identifiability}: We prove that
$\eta$, $q$, and $\Ast$ are uniquely determined
by the observations under a rank condition on
the increments of the observed attractiveness
from four selected trajectory--time pairs. ii) \emph{Information loss from an unknown $\Ast$}: 
  {We quantify the information about $\eta$ that is lost when $\Ast$ is unknown and must be inferred jointly.}
  iii) \emph{Information loss near a no-flux boundary}: We show that, near a smooth no-flux boundary, the information about $\eta$ decreases at explicit rates, both in the data and in the discrete learner. iv) \emph{Numerical verification}: On a square and on the geometry of Chicago, the proposed method recovers all three parameters with relative errors of a few percent, and the predicted losses of information are observed.

%
%


Next, we relate our approach to existing work.
Two classical approaches to estimate {unknown parameters in PDEs} are the
output least-squares approach and the equation-error approach.  In the first,
the misfit between the data and the solution of the forward problem is
minimized over the parameters, possibly with regularization or in a Bayesian
setting~\cite{banks1989estimation,stuart2010bayesian}. Each evaluation
requires a forward solve, and the problem is in general nonconvex.  In the
second, the observed fields are substituted into the PDE and the residual is
minimized over the parameters~\cite{richter1981inverse,falk1983error,kohn1988variational}.
Our method
belongs to this second class, with the residual measured in a weak norm so
that only first derivatives of the data are needed.

The use of weak formulations in parameter estimation is closely related to their development in data-driven identification of differential equations, where sparse regression over a library of candidate terms~\cite{brunton2016sindy,rudy2017pdefind} was combined with weak formulations to avoid pointwise derivatives of noisy data~\cite{gurevich2019robust,reinbold2020noisy,messenger2021wsindypde,messenger2021wsindymms,tang2023weakident}. When the model structure is known, weak-form parameter estimation, as in WENDy~\cite{bortz2023wendy}, is often faster and more accurate than forward-solver-based least squares. Spatially or temporally varying coefficients have been treated in~\cite{rudy2019parametric,kang2021ident} and, closest to our setting, in~\cite{tang2025wgident}, where the coefficients and the test functions are represented with B-splines. 

The choice of test functions plays an important role in weak-form learning: they have been selected using spectral information~\cite{messenger2021wsindypde}, at multiple scales~\cite{bortz2023wendy}, based on quadrature errors~\cite{tran2026weak}, or made dependent on the unknown parameters~\cite{gao2026self}. In our method, the spatial test functions are chosen from the space used to represent the parameters, which simplifies the assembly of the discrete dual norm and provides a compatibility relation for the unknown static attractiveness, although unique solvability of the learner still requires sufficient test-space resolution.


Our analysis is close in spirit to learning interaction kernels of agent-based systems from trajectory data, where a coercivity condition guarantees identifiability and optimal convergence rates~\cite{mauro2017,lu2019pnas, feng2024learning}, with related estimators for stochastic differential equations~\cite{guo2024sde,chen2026kellersegel} and mean-field equations~\cite{lang2022meanfield,messenger2022meanfield}. Our model is itself the mean-field limit of an agent-based model~\cite{short2008statistical,Hao2026}, and the information coercivity of our learner (see Sec.~\ref{sec:boundary-information}) plays the role of that coercivity condition.

Statistical models of crime, such as self-exciting point processes, are
routinely fitted to event data by maximum
likelihood~\cite{mohler2020impact,park2021investigating} or sequential data
assimilation~\cite{santitissadeekorn2018sequential}.  For the PDE model
of~\cite{short2008statistical}, 
constant parameters and initial conditions
were fitted to simulated crime data through a
likelihood based on the predicted crime rate in~\cite{lloyd2016exploring}, and this likelihood was found to be
flat over large ranges of the parameters.  This is an early indication that
identifiability is a central issue for this model.  To our knowledge, the
recovery of spatially varying $\eta$, $\Ast$, and $q$, and the question of
when they are identifiable, have not been addressed before.

This paper is {Part~I} of a series.  Here, the numerical results use
idealized observations: noise-free, densely sampled fields of both the
attractiveness and the criminal density.  Part~II will extend the model to
include police control, following \cite{Hao2026b}, and study the sensitivity
 of the method to observation noise and {resolution}. {Part~III will address the inverse problem with indirect observations, where only the crime rate is observed rather than the attractiveness and criminal density separately.}

The rest of the paper is organized as follows. Sec.~\ref{sec:pd} recalls the PDE model and specifies the observations. Sec.~\ref{sec:weak-form} introduces the weak-form learning method and gives conditions for the solvability of the learner. Sec.~\ref{sec:ident} establishes identifiability and quantifies the information lost when $\Ast$ is unknown.  Sec.~\ref{sec:boundary-information} analyzes the loss of information near a no-flux boundary. Numerical results are presented in Sec.~\ref{sec:num_res}, and conclusions are drawn in Sec.~\ref{sec:concl}.
%
\section{Problem description}\label{sec:pd}
This section briefly recalls the model from \cite{short2008statistical} written in the formulation we have used in \cite{Hao2026}. 
We consider a bounded, connected Lipschitz domain $\Omega\subset \mathbb{R}^2$ and a time interval of interest $(0,T)$.
Let $A(t, \bx)>0$ denote the level of attractiveness after suitable nondimensionalization. This variable, which measures the burglars’
perception at time $t$ of the attractiveness of a site located at $\bx$, is the sum of a static component $\Ast(\bx)$ and a dynamic component associated with repeat
and near-repeat victimization. Let 
$\rho(t, \bx)$ be the density
of criminals after suitable nondimensionalization. Then, the non-dimensional PDE model from \cite{short2008statistical} can be written as: find $A(t, \bx)$ and $\rho(t, \bx)$ such that 
\begin{align}
    &\partial_t A(t, \bx) - \nabla\cdot(\eta(\bx)\nabla A(t, \bx)) + A(t, \bx) - \rho(t, \bx)A(t, \bx) = -\nabla\cdot(\eta(\bx)\nabla \Ast(\bx)) + \Ast(\bx), \label{eq:PDE1} \\
    &\partial_t \rho(t, \bx) - \nabla\cdot(\nabla\rho(t, \bx) - 2\rho(t, \bx)\nabla(\log(A(t, \bx)))) + \rho(t, \bx)A(t, \bx) = q(\bx), \label{eq:PDE2}
\end{align}
where {$\partial_t$ denotes differentiation with respect to time},
$\eta(\bx)$ characterizes the strength of the 
near-repeat victimization and the broken-window effect and $q(\bx)$ is a source term related to criminal regeneration and time scale for repeat victimization. We assume that
\begin{equation}
 0<\eta_{\min}\leq\eta(\bx)\leq\eta_{\max}<\infty.
 \label{eq:eta-positive}
\end{equation}

Problem \eqref{eq:PDE1}-\eqref{eq:PDE2}
is supplemented with boundary conditions 
\begin{equation}\label{eq:BC}
    \partial_{\mathbf n}\rho = \nabla\rho(t, \bx)\cdot\mathbf{n} = 0, \quad \partial_{\mathbf n} A = \nabla A(t, \bx)\cdot\mathbf{n} = 0, \quad \text{on } \partial\Omega,
\end{equation}
and initial conditions 
\[
\rho(0, \bx) = \rho_0(\bx), \quad A(0, \bx) = A_0(\bx).
\]
In addition, we assume that 
\begin{equation}\label{eq:nofluxAst}
  \nabla \Ast(\bx)\cdot\mathbf{n} = 0, \quad \text{on } \partial\Omega.
\end{equation}

While model \eqref{eq:PDE1}-\eqref{eq:PDE2} is 
versatile and able to reproduce patterns observed in real crime data, quantities like
$\eta(\bx)$, $q(\bx)$, and $\Ast(\bx)$
are hard to measure in practice. Hence, we will assume that they are unknown.
We consider two cases throughout the paper: $\Ast$ known, in which case only $\eta$ and $q$ are learned, and $\Ast$ unknown, in which case three parameter fields are learned from only two equations.

To be able to learn functions $\eta(\bx)$, $q(\bx)$, and $\Ast(\bx)$, {we assume that direct observations of $A$ and $\rho$
are available.} More specifically, let  $i=1,\ldots,N$, denote a trajectory, where $N$ is the total number of trajectories.  By ``trajectory'', we mean an evolution of the system from a given initial condition, so different trajectories are associated with different initial conditions.  For each trajectory $i$, the observed fields are the attractiveness and criminal density, denoted with $A_{(i)}(t,\bx)$ and $\rho_{(i)}(t,\bx)$. Clearly, 
$A_{(i)}(t,\bx)$ and $\rho_{(i)}(t,\bx)$ need to satisfy 
\eqref{eq:PDE1}-\eqref{eq:PDE2} and 
\eqref{eq:BC}. Since we will use a learning method for the model parameters based on the weak formulation of problem \eqref{eq:PDE1}-\eqref{eq:PDE2}, it is sufficient to require that $A_{(i)}$ is bounded below by a positive
constant and that
\begin{equation}
 A_{(i)},\rho_{(i)}\in L^2(0,T;H^1(\Om)),
 \qquad
 A_{(i)}^2\nabla(\rho_{(i)}/A_{(i)}^2)\in L^2((0,T)\times\Om)^2.
 \label{eq:weak-data-regularity}
\end{equation}
For the same reason, integrability of the source terms in \eqref{eq:PDE1}-\eqref{eq:PDE2} is required too.
Higher regularity of the observed fields is assumed only in the analysis: in
Sec.~\ref{sec:ident}, for the identifiability result of
Th.~\ref{thm:identifiability} and its consequences, and in
Sec.~\ref{sec:boundary-information}, for the behavior of the information
near the boundary.  It is not needed to formulate or apply the learning
method.

In practice, $A$ and $\rho$ are not observed directly: crime data count events, whose expected rate is related to $\rho A$~\cite{Hao2026b}. Here, the observations are synthetic (see Sec.~\ref{sec:num_res}), while indirect observations will be addressed in Part~III.

\section{A weak-form learning method}\label{sec:weak-form}
When $\Ast$ is known, equations \eqref{eq:PDE1}--\eqref{eq:PDE2} are linear
in the unknown parameters: \eqref{eq:PDE1} in $\eta$ and \eqref{eq:PDE2} in
$q$.  When $\Ast$ is unknown, \eqref{eq:PDE1} is no longer linear.  To
recover linearity, we collect the time-independent terms into the
{static forcing}
\begin{equation}
 S:=\Leta\Ast,
 \qquad\text{where}\quad
 \Leta v:=v-\diver(\eta\nabla v),
 \label{eq:static-source}
\end{equation}
so that \eqref{eq:PDE1} becomes
\begin{equation}
 \partial_tA-\diver(\eta\nabla A)+A-\rho A=S,
 \label{eq:static-source-A}
\end{equation}
which is jointly linear in $(\eta,S)$. 

Let $\phi(t,\bx)$ be a space--time test function in $L^2(0,T;H^1(\Om))$, with
$\partial_t\phi$ square integrable, that vanishes at $t=0$ and $t=T$.
For each trajectory $i=1,\ldots,N$, the weak forms of
\eqref{eq:PDE2} and \eqref{eq:static-source-A} tested against $\phi$ are 
\begin{align}
 &\int_0^T\!\int_\Om
 \left(\eta\nabla A_{(i)}\cdot\nabla\phi-S\phi\right)\dd \bx\dd t
 =
 \int_0^T\!\int_\Om
 \left(A_{(i)}\partial_t\phi-(A_{(i)}-\rho_{(i)}A_{(i)})\phi\right)\dd \bx\dd t,
 \label{eq:weak-joint}\\
 &\int_0^T\!\int_\Om (q\phi)\dd \bx\dd t
 =
 \int_0^T\!\int_\Om
 \left( -\rho_{(i)}\partial_t\phi
 +A_{(i)}^2\nabla\!\left(\frac{\rho_{(i)}}{A_{(i)}^2}\right)\cdot\nabla\phi
 +\rho_{(i)}A_{(i)}\phi \right) \dd \bx\dd t.
 \label{eq:weak-q}
\end{align}
Because of the no-flux conditions \eqref{eq:BC}, integrating by parts in
space produces no boundary terms, so the weak forms
\eqref{eq:weak-joint}--\eqref{eq:weak-q} involve neither the values of the
observed fields nor their derivatives on $\partial\Om$. 

Let $q\in L^2(\Omega)$ and $\eta\in L^\infty(\Omega)$ with \eqref{eq:eta-positive}, and recall that $S\in H^{-1}(\Omega)$.  {For simplicity of notation}, the same parameter space is used for all three parameters and we choose it conforming for all of them:
$X_J=\operatorname{span}\{\psi_j\}_{j=1}^J\subset L^\infty(\Omega)$. 
We write the model parameters as a linear combination of functions $\psi_j (\bx)$, $j=1, \dots, J$: 
\begin{equation}\label{eq:FD_exp}
   \eta_J=\sum_{j=1}^J\alpha_j\psi_j,
 \qquad
 S_J=\sum_{j=1}^Js_j\psi_j,
 \qquad
 q_J=\sum_{j=1}^J\beta_j\psi_j.  
\end{equation}
{Of course, one could use three different parameter spaces. The choice of using the same does not affect the analysis or computations.}
Let $V_K=\operatorname{span}\{\varphi_\ell\}_{\ell=1}^K$ be the space--time test space. Substituting \eqref{eq:FD_exp} into \eqref{eq:weak-joint}--\eqref{eq:weak-q} and testing each trajectory $i=1,\ldots,N$ with $\varphi_\ell$, $\ell=1,\ldots,K$, gives one equation for every pair $(i,\ell)$ and hence $NK$ equations. Next, we distinguish the two settings of 
$\Ast$ known and unknown.

If $\Ast$ is known, we learn $(\eta,q)$ by solving the following matrix problems:
\begin{align}
&M_a\boldsymbol\alpha=\bb_a,\label{eq:MKprob} \\
&M_q\boldsymbol\beta=\bb_\rho,
 \label{eq:Mqprob}
\end{align}
where $M_a,M_q\in\mathbb R^{NK\times J}$, $\bb_a,\bb_\rho\in\mathbb R^{NK}$, and $\boldsymbol\alpha=(\alpha_j)_{j=1}^J$ and $\boldsymbol\beta=(\beta_j)_{j=1}^J$ are the vectors containing the unknown coefficients. The rows in the matrix problems come from stacking distinct trajectories.  The matrices are generally rectangular, thus problems  \eqref{eq:MKprob} and \eqref{eq:Mqprob} have to be solved in the least-squares sense.
With rows indexed by $\tilde\ell$ with $\tilde\ell = (i - 1)K + \ell$ and columns indexed by $j$, the entries of the matrices are 
\begin{align}
 (M_a)_{\tilde\ell,j}
 &=\frac1{NT}\int_0^T\!\int_\Om
   \psi_j\nabla(A_{(i)}-\Ast)\cdot\nabla\varphi_\ell\,\dd\bx\dd t,\label{eq:MK}\\
 (M_q)_{\tilde\ell,j}
 &=\frac1{NT}\int_0^T\!\int_\Om
   \psi_j\varphi_\ell\,\dd\bx\dd t, \label{eq:Mq}
\end{align}
and the entries of the right-hand side vectors are
\begin{align*}
 (\bb_a)_{\tilde\ell}
 &=\frac1{NT}\int_0^T\!\int_\Om
   \left[A_{(i)}\partial_t\varphi_\ell
   -(A_{(i)}-\rho_{(i)}A_{(i)})\varphi_\ell\right]\dd\bx\dd t+\frac1{NT}\int_0^T\!\int_\Om
   \Ast\varphi_\ell\,\dd\bx\dd t,\\
 (\bb_\rho)_{\tilde\ell}
 &=\frac1{NT}\int_0^T\!\int_\Om
   \left[-\rho_{(i)}\partial_t\varphi_\ell
   +A_{(i)}^2\nabla\!\left(\frac{\rho_{(i)}}{A_{(i)}^2}\right)
     \cdot\nabla\varphi_\ell
   +\rho_{(i)}A_{(i)}\varphi_\ell\right]\dd\bx\dd t.
\end{align*}

If $\Ast$ is unknown, we learn $(\eta,S,q)$ by solving the following matrix problem:
\begin{equation}\label{eq:MeMs}
[M_\eta\;M_S]
\begin{pmatrix}\boldsymbol\alpha\\\boldsymbol s\end{pmatrix}=\bb_\eta,
\end{equation}
and \eqref{eq:Mqprob}. Here, $M_\eta,M_S\in\mathbb R^{NK\times J}$,
$[M_\eta\;M_S]\in\mathbb R^{NK\times2J}$,
$\bb_\eta\in\mathbb R^{NK}$, and $\boldsymbol s=(s_j)_{j=1}^J$ is the remaining  coefficient vector from 
\eqref{eq:FD_exp}.
With the same row index $\tilde\ell$ and column index $j$, the entries for the matrices are 
\begin{align}
 (M_\eta)_{\tilde\ell,j}
 &=\frac1{NT}\int_0^T\!\int_\Om
   \psi_j\nabla A_{(i)}\cdot\nabla\varphi_\ell\,\dd\bx\dd t, \label{eq:Meta} \\
 (M_S)_{\tilde\ell,j}
 &=-\frac1{NT}\int_0^T\!\int_\Om
   \psi_j\varphi_\ell\,\dd\bx\dd t, \label{eq:Ms}
\end{align}
and the entries of the right-hand side vector are
\begin{align*}
 (\bb_\eta)_{\tilde\ell}
 &=\frac1{NT}\int_0^T\!\int_\Om
   \left[A_{(i)}\partial_t\varphi_\ell
   -(A_{(i)}-\rho_{(i)}A_{(i)})\varphi_\ell\right]\dd\bx\dd t.
\end{align*}
Then, $\Ast$ is easily recovered from $S=\Leta\Ast$. {More details on this are presented below, before Rem.~\ref{rem:pos_eta}.}

The residuals of the matrix problems above should be measured in a norm that reflects the functions they come from, rather than in the Euclidean norm, which depends on the choice of test functions.  Since the test functions belong to $L^2(0,T;H^1(\Om))$, we use the corresponding dual norm. This norm is built from the following Gram matrix of the test functions:
\begin{align}
 \mathsf Q_{\tilde\ell, \tilde\ell'}
 &=\frac{\delta_{ii'}}{NT}\int_0^T\!\int_\Om
 \left(\varphi_\ell\varphi_{\ell'}
+\nabla\varphi_\ell\cdot\nabla\varphi_{\ell'}\right)
 \,\dd \bx\dd t,
 \label{eq:test-Gram}
\end{align}
where $\delta_{ii'}$ is the Kronecker delta {and $\tilde\ell' = (i' - 1)K + \ell'$}.  
The matrix $\mathsf Q\in\R^{NK\times NK}$ is block diagonal across trajectories, since it uses the same stacking strategy as matrices \eqref{eq:MK}, \eqref{eq:Mq}, \eqref{eq:Meta}, and \eqref{eq:Ms}. It is also normalized like those matrices, i.e., every entry is divided by $NT$. For a residual vector $\br\in\R^{NK}$, this norm is given by $\norm{\br}_{\mathsf Q^{-1}}^2=\br^\top\mathsf Q^{-1}\br$, a discrete version of the squared norm in $L^2(0,T;H^{-1}(\Om))$.

We solve problems \eqref{eq:MKprob}, \eqref{eq:Mqprob}, and \eqref{eq:MeMs} in the least-squares sense using this discrete norm.
When $\Ast$ is known, this means solving
\begin{align}
 \argmin_{\boldsymbol{\alpha}}
 \norm{M_a\boldsymbol{\alpha}-\bb_a}_{\mathsf Q^{-1}}^2
 \quad \text{and} \quad
 \argmin_{\boldsymbol{\beta}}
 \norm{M_q\boldsymbol{\beta}-\bb_\rho}_{\mathsf Q^{-1}}^2.
 \label{eq:weighted-learner_AstK}
\end{align}
When $\Ast$ is unknown, this means solving
\begin{align}
 \argmin_{\boldsymbol{\alpha},\boldsymbol{s}}
 \norm{M_\eta\boldsymbol{\alpha}+M_S \boldsymbol{s}-\bb_\eta}_{\mathsf Q^{-1}}^2
 \quad \text{and} \quad
 \argmin_{\boldsymbol{\beta}}
 \norm{M_q\boldsymbol{\beta}-\bb_\rho}_{\mathsf Q^{-1}}^2.
 \label{eq:weighted-known-q-learner_AstU}
\end{align}
We refer to the least-squares problems \eqref{eq:weighted-learner_AstK} and
\eqref{eq:weighted-known-q-learner_AstU} as the {learner}. More
specifically, the {diffusivity learner} is 
the problem for $\eta$ (together with $S$ when $\Ast$ is unknown) 
and the {source learner}
is the problem for
$q$.

After solving \eqref{eq:weighted-known-q-learner_AstU},
we obtain $(\eta_J,S_J)$ as defined in \eqref{eq:FD_exp}. The static attractiveness $\Ast_J$ is then recovered by solving the following  elliptic problem: find $\Ast_J\in H^1(\Om)$ such that
\begin{equation}
 \int_\Om \Ast_Jv+\eta_J\nabla\Ast_J\cdot\nabla v\,\dd \bx
 =\ip{S_J}{v},
 \qquad\forall v\in H^1(\Om).
 \label{eq:recover-Ast}
\end{equation}

\begin{remark}\label{rem:pos_eta}
The parameter space $X_J$ does not enforce positivity
of $\eta_J$ or the recovered $\Ast_J$.
Positivity of $\eta_J$ can be imposed through coefficient
constraints, such as $\alpha_j\geq\eta_{\min}$ for
B-spline bases.
We solve the unconstrained problems since the learned
$\eta_J$ satisfies \eqref{eq:eta-positive} in all
reported experiments. This also guarantees that \eqref{eq:recover-Ast} is coercive.
\end{remark}
If $\Ast$ is unknown, $S$ is learned together with $\eta$.  For a fixed
coefficient vector $\boldsymbol\alpha$, the residual
$M_\eta\boldsymbol\alpha+M_S\boldsymbol s-\bb_\eta$ is linear in
$\boldsymbol s$, so the optimal $\boldsymbol s$ is given by a linear
least-squares fit: it cancels the part of $M_\eta\boldsymbol\alpha-\bb_\eta$
that can be written as $M_S\boldsymbol s$, i.e., the part that can be
matched by a suitable choice of the static forcing.  Only the remaining part of the residual constrains
$\eta$.  The cancelled part is given by the $\mathsf Q^{-1}$-orthogonal
projector onto the range of $M_S$:
\begin{equation}
 \mathsf P_S
 :=M_S(M_S^\top \mathsf Q^{-1} M_S)^{\dagger} M_S^\top \mathsf Q^{-1},
 \label{eq:discrete-static-projector}
\end{equation}
where ${}^{\dagger}$ denotes the Moore--Penrose pseudoinverse.  

\begin{proposition}
\label{prop:exact-discrete-solvability}
Suppose that $\mathsf Q$ is positive definite, which holds whenever the
test functions are linearly independent.
If $\Ast$ is known, the diffusivity learner is unique if and only if
\begin{equation}
 M_a^\top \mathsf Q^{-1}M_a\succ0.
 \label{eq:HK-discrete}
\end{equation}
If $\Ast$ is unknown, the diffusivity learner is unique if and
only if
\begin{equation}
 M_\eta^\top \mathsf Q^{-1}(I-\mathsf P_S)M_\eta\succ0.
 \label{eq:HU-discrete}
\end{equation}
In both cases of known and unknown $\Ast$, the source learner is unique if and only if
\begin{equation}
 M_q^\top \mathsf Q^{-1}M_q\succ0.
 \label{eq:Hq-discrete}
\end{equation}
\end{proposition}

\begin{proof}
The weighted least-squares problem
$\norm{A\bc-\bb}_{\mathsf Q^{-1}}^2$
has a unique minimizer if and only if
$A^\top\mathsf Q^{-1}A\succ0$.
Taking $A=M_a$ and $A=M_q$ gives
\eqref{eq:HK-discrete} and \eqref{eq:Hq-discrete},
respectively.
If $\Ast$ is unknown, minimizing over $\boldsymbol s$
gives $A=(I-\mathsf P_S)M_\eta$.
The $\mathsf Q^{-1}$-orthogonality of $\mathsf P_S$
then yields \eqref{eq:HU-discrete}.
\end{proof}

{
In our numerical experiments and parts of the analysis in
Sec.~\ref{sec:ident}, we use {spatial self-testing},
where the spatial test functions are chosen from the
same basis used to represent the unknown parameters.
We use
$\varphi_{(m,p)}(t,\bx)=\theta_m(t)\psi_p(\bx)$,
where $\theta_m\in H^1_0(0,T)$,
$m=1,\ldots,M_t$, are temporal test functions and
$\psi_p$, $p=1,\ldots,J$, are the spatial basis functions.
Each trajectory then has $K = M_tJ$ test functions.
The temporal and spatial mass matrices and the
spatial stiffness matrix are
\begin{equation}
 (\mathsf Q_t)_{mn}=\int_0^T\theta_m\theta_n\,\dd t,
 \qquad
 \mathsf R_{pq}=\int_\Om\psi_p\psi_q\,\dd\bx,
 \qquad
 (\mathsf K_x)_{pq}=\int_\Om\nabla\psi_p\cdot\nabla\psi_q\,\dd\bx,
 \label{eq:mass-stiffness}
\end{equation}
so the diagonal blocks of $\mathsf Q$ in
\eqref{eq:test-Gram} are
\begin{equation}
 \frac1{NT}\,\mathsf Q_t\otimes(\mathsf R+\mathsf K_x).
 \label{eq:self-test-grams}
\end{equation}
This choice also gives the discrete compatibility
relation established in
Lem.~\ref{lem:self-test-compatibility}.
}

\section{Identifiability and the effect of an unknown $\Ast$}\label{sec:ident}

Prop.~\ref{prop:exact-discrete-solvability} gives rank
conditions for the learner.  In this section, we derive a sufficient
condition for identifiability at the PDE level, namely a rank condition on the
observed transients (Th.~\ref{thm:identifiability}).
Here, identifiability means that the observations determine the parameters
uniquely: if two admissible triples $(\eta_j,q_j,\Ast_j)$, $j=1,2$, generate
the same observed fields $(A_{(i)},\rho_{(i)})_{i=1}^N$ through
\eqref{eq:PDE1}--\eqref{eq:PDE2}, then $\eta_1=\eta_2$, $q_1=q_2$, and
$\Ast_1=\Ast_2$ almost everywhere.
Then, we compare the
information available when $\Ast$ is known or unknown, and show how this
PDE-level rank condition implies the discrete conditions of
Prop.~\ref{prop:exact-discrete-solvability}.

\subsection{Continuum identifiability 
}

We consider the harder setting in which $\Ast$ is unknown. {The uniqueness result below also holds when $\Ast$ is known too.}

Recall the elliptic operator $\Leta$ in \eqref{eq:static-source}, which we consider
with homogeneous Neumann boundary condition.  
In its variational form $\langle  \Leta v, u \rangle$, this operator is associated with the bilinear form
\[
 a_\eta(u,v)=\int_\Om uv+\eta\nabla u\cdot\nabla v\,\dd \bx, \quad u,v \in H^1(\Om).
\]
We also introduce the dual space
\(
 H^{-1}(\Om):=(H^1(\Om))',
\)
because weak residuals act on $H^1(\Om)$ test functions without a homogeneous
boundary trace. By \eqref{eq:eta-positive}, $a_\eta$ is coercive on $H^1(\Om)$, so that, in the distributional sense,
\begin{equation}
 \Leta:H^1(\Om)\longrightarrow H^{-1}(\Om)
 \quad\text{is an isomorphism.}
 \label{eq:Leta-isomorphism}
\end{equation}


Choose four distinct trajectory--time pairs
$(i_k,t_k)$, $k=0,1,2,3$, at which
\eqref{eq:PDE1} holds.
Define three increments
\begin{equation}\label{eq:increments}
\delta A_k(\bx)
:=A_{(i_k)}(t_k,\bx)-A_{(i_0)}(t_0,\bx),
\qquad k=1,2,3,
\end{equation}
and the local increment matrix $\mathsf G(\bx)\in\R^{3\times 3}$ as follows
\begin{equation}
 \mathsf G(\bx):=
 \begin{pmatrix}
  \partial_{x_1}(\delta A_1)&\partial_{x_2}(\delta A_1)&\Delta (\delta A_1)\\
  \partial_{x_1}(\delta A_2)&\partial_{x_2}(\delta A_2)&\Delta (\delta A_2)\\
  \partial_{x_1}(\delta A_3)&\partial_{x_2}(\delta A_3)&\Delta (\delta A_3)
 \end{pmatrix}.
 \label{eq:rank3matrix}
\end{equation} 

Let us now suppose that two admissible triples $(\eta_j,q_j,\Ast_j)$, $j=1,2$, generate the same observations, and let $u:=\eta_1-\eta_2$. In the proof of Th.~\ref{thm:identifiability}, we will show that the terms involving $\Ast_1$ and $\Ast_2$ are independent of the trajectory and time, and therefore cancel in the increments \eqref{eq:increments}, so that $u$ satisfies
\[
 \diver(u\nabla\delta A_{k})=0, \qquad k=1,2,3.
\]
If $u$ were differentiable, expanding these equations would give
$
 \mathsf G(\bx)\,[\partial_{x_1}u,\ \partial_{x_2}u,\ u]^\top=\boldsymbol{0},
$
and if $\mathsf G(\bx)$ has rank $3$ (i.e., it has full column-rank) for almost every $\bx\in\Om$, this would force $u=0$.  Since the diffusivities are only assumed to be bounded, we cannot differentiate $u$.  The following lemma reaches the same conclusion in weak form. 

\begin{lemma}\label{lem:weak-cofactor}
Let $f_1,f_2,f_3\in W^{2,\infty}_{\rm loc}(\Om)$, and let $\mathsf G$ be the
matrix \eqref{eq:rank3matrix} with $\delta A_k$ replaced by $f_k$.  Let
$\mathbf w=(w_k)_{k=1}^3$ be the cross product of the first two columns of
$\mathsf G$.  Then, for every $u\in L^2_{\rm loc}(\Om)$,
\begin{equation}
 \int_\Om u\det\mathsf G\,\phi\,\dd\bx=\sum_{k=1}^3\langle \diver(u\nabla f_k),w_k\phi\rangle
 \qquad\forall\phi\in C_c^\infty(\Om).
 \label{eq:weak-cofactor-identity}
\end{equation}
In particular, if $\diver(u\nabla f_k)=0$ for $k=1,2,3$ and $\det\mathsf G\ne0$ almost
everywhere, then $u=0$ almost everywhere.
\end{lemma}

\begin{proof}
The entries of $\mathbf w$ are products of first derivatives of the $f_k$, so
$w_k\in W^{1,\infty}_{\rm loc}(\Om)$.  Since $\mathbf w$ is orthogonal to the
first two columns of $\mathsf G$, we have $\sum_kw_k\nabla f_k=0$.  Taking its
divergence and expanding $\det\mathsf G$ along the third column gives
$\sum_k\nabla w_k\cdot\nabla f_k=-\sum_kw_k\Delta f_k=-\det\mathsf G$
almost everywhere.  For any $\phi\in C_c^\infty(\Om)$, the products
$w_k\phi$ are admissible $H_0^1$ test functions on a compact subdomain.
Consequently,
\[
 \sum_k\langle \diver(u\nabla f_k),w_k\phi\rangle
 =-\int_\Om u\sum_k\nabla f_k\cdot\nabla(w_k\phi)\,\dd\bx
 =\int_\Om u\det\mathsf G\,\phi\,\dd\bx.
\]
This proves the identity without differentiating $u$.
\end{proof}

\begin{theorem}
\label{thm:identifiability}
Let two admissible triples $(\eta_j,q_j,\Ast_j)$, $j=1,2$, generate
the same observations $(A_{(i)},\rho_{(i)})$.  Assume
$\eta_j\in L^\infty(\Om)$ satisfy \eqref{eq:eta-positive},
$\Ast_j\in H^1(\Om)$, the observed fields satisfy \eqref{eq:PDE1}--\eqref{eq:PDE2} and
\eqref{eq:BC}, and the observed attractiveness fields
$A_{(i_k)}(t_k,\cdot)$ at the four selected trajectory--time pairs, $k=0,1,2,3$,
belong to $W^{2,\infty}_{\rm loc}(\Om)$.
If
\begin{equation}
 \operatorname{rank}\mathsf G(\bx)=3
 \quad\text{for almost every }\bx\in\Om,
 \label{eq:rank3condition}
\end{equation}
then $\eta_1=\eta_2$, $q_1=q_2$, and $\Ast_1=\Ast_2$ almost everywhere.
\end{theorem}

\begin{proof}
For every trajectory $i$, the common observations determine
$q=\partial_t\rho_{(i)}-\diver\bigl(\nabla\rho_{(i)}-2\rho_{(i)}\nabla\log A_{(i)}\bigr)+\rho_{(i)}A_{(i)}$
as a distribution, so $q_1=q_2$.  Subtracting the attractiveness equations
with $u=\eta_1-\eta_2$ gives
\begin{equation}
 -\diver\bigl(u\nabla A_{(i_k)}(t_k,\cdot)\bigr)
=\mathcal L_{\eta_1}\Ast_1-\mathcal L_{\eta_2}\Ast_2,
\qquad k=0,1,2,3.
 \label{eq:eta-difference-static}
\end{equation}
The right-hand side depends only on the parameters, which are shared by
all trajectories, so it is independent of the trajectory and time.
Taking three differences yields
\begin{equation}
 \diver(u\nabla\delta A_k)=0,\qquad k=1,2,3.
 \label{eq:increment-equations}
\end{equation}
Lem.~\ref{lem:weak-cofactor} and \eqref{eq:rank3condition} imply $u=0$.
The remaining difference is $\Leta(\Ast_1-\Ast_2)=0$.  Testing with
$\Ast_1-\Ast_2$, and using coercivity, gives $\Ast_1=\Ast_2$.
\end{proof}

\begin{remark}\label{rem:multi-snapshot-rank}
The four trajectory--time pairs may depend on $\bx$. If $\mathsf G$ has rank
three only on a subset $\Om'\subset\Om$, \eqref{eq:weak-cofactor-identity}
still gives $u=0$ almost everywhere on $\Om'$. Hence, the conclusion of
Th.~\ref{thm:identifiability} holds whenever $\Om$ is covered, up to a set
of measure zero, by countably many such subsets, each with its own choice of
four pairs satisfying the regularity assumption of the theorem.
\end{remark}

The identifiability of a diffusion coefficient, such as $\eta$, from interior data is a
classical question.  With a single observed solution of a steady diffusion
equation, the coefficient is uniquely determined, provided that it is
prescribed on the inflow part of the boundary and that the gradient and the
Laplacian of the solution do not degenerate simultaneously~\cite{richter1981inverse}.
In hybrid inverse problems, coefficients are reconstructed from several
interior solutions whose gradients span all directions at every
point~\cite{bal2013reconstruction}.  The rank condition
\eqref{eq:rank3condition} combines both ideas, since it involves the
gradients and the Laplacians of the increments \eqref{eq:increments} of
several observed fields.  
However, unlike the results in \cite{richter1981inverse,bal2013reconstruction},
Th.~\ref{thm:identifiability} requires no prescribed values of $\eta$,
$\Ast$, and $q$ anywhere inside or on the boundary of the domain.

\subsection{Information loss from an unknown $\Ast$}
Th.~\ref{thm:identifiability} shows that, under the rank condition
\eqref{eq:rank3condition}, the diffusivity $\eta$ is identifiable even when $\Ast$
is unknown.  However,
it does not quantify how much information is lost by
having to learn $\Ast$.  The goal of this section is to quantify the loss.

When $\Ast$ is unknown, the static forcing
$S=\Leta\Ast$ in \eqref{eq:static-source} is effectively a free unknown:
since $\Leta$ is invertible by \eqref{eq:Leta-isomorphism}, every
$S\in H^{-1}(\Om)$ arises from some $\Ast$, for any $\eta$.  Hence, any
component of the residual that does not vary with the trajectory--time pair $(i,t)$ can be matched by a
suitable choice of $S$, whatever $\eta$ is, and therefore carries no
information about $\eta$.
We want to compare how strongly the data constrain $\eta$
when $\Ast$ is known and when it is unknown. For this purpose, we measure the residual of \eqref{eq:PDE1}, i.e.,
the amount $F_{(i)} (t, \cdot)$ by which the attractiveness equation fails for candidate
parameters at each $(i,t)$ pair.  Since this residual
contains divergence terms, it lies in $H^{-1}(\Om)$, and we use the space
\[
 Y:=L^2\!\left(0,T;H^{-1}(\Om)\right)^N,
 \qquad
 \norm{F}_Y^2
 =\tavg{\norm{F_{(i)} (t,\cdot)}_{H^{-1}}^2}.
\]

Because the residual is affine in $\eta$, perturbing $\eta$ by $u\in L^\infty(\Om)$ changes it
by the following map 
\begin{equation}
 (\Tcal_Ku)_{(i)}(t,\cdot):=-\diver\!\left(u\nabla\bigl(A_{(i)}(t,\cdot)-\Ast\bigr)\right),
 \label{eq:TK}
\end{equation}
in the case of known $\Ast$. The subscript $K$ in \eqref{eq:TK} underscores the fact that $\Ast$ is known. Below, we will use subscript $U$ for the case of unknown $\Ast$.
Without knowledge of $\Ast$, the residual changes instead by
$-\diver(u\nabla A_{(i)}(t,\cdot))$, which differs from $\Tcal_Ku$ by the term
$-\diver(u\nabla\Ast)$.  This term is independent of trajectory--time pair $(i,t)$, so it can be
absorbed by the unknown static forcing $S$.




\begin{theorem}
\label{thm:profiling}
Let
\[
 \bar A:=\tavg{A_{(i)}(t,\cdot)}
\]
be the average of the observed attractiveness.  For every $u\in L^\infty(\Om)$,
\begin{equation}
 \inf_{s\in H^{-1}(\Om)}\norm{\Tcal_Ku-s}_Y=\norm{\Tcal_Uu}_Y,
 \qquad
 (\Tcal_Uu)_{(i)}(t,\cdot):=-\diver\!\left(u\nabla\bigl(A_{(i)}(t,\cdot)-\bar A\bigr)\right),
 \label{eq:TU-centered}
\end{equation}
and
\begin{equation}
 \norm{\Tcal_Ku}_Y^2
 =\norm{\Tcal_Uu}_Y^2
 +\norm{\diver\!\left(u\nabla(\bar A-\Ast)\right)}_{H^{-1}}^2.
 \label{eq:information-ordering}
\end{equation}
In particular, $\norm{\Tcal_Uu}_Y\leq\norm{\Tcal_Ku}_Y$, i.e., {information on $\eta$ may be lost when $\Ast$ is unknown}.  What is lost is exactly the
response to the averaged field $\bar A-\Ast$. Only the variation of the
observations around their average informs $\eta$.
\end{theorem}

\begin{proof}
The minimizer of $\norm{\Tcal_Ku-s}_Y$ over
$s\in H^{-1}(\Om)$ is the trajectory--time average
\[
 s=-\diver\!\left(u\nabla(\bar A-\Ast)\right).
\]
Subtracting this average gives $\Tcal_Ku-s=\Tcal_Uu$,
which proves \eqref{eq:TU-centered}.
Since $\Tcal_Uu$ has zero trajectory--time average,
it is orthogonal to $s$ in $Y$.
This gives \eqref{eq:information-ordering}.
\end{proof}





Th.~\ref{thm:profiling} concerns an unrestricted static forcing
$s\in H^{-1}(\Om)$.  In the learner of Sec.~\ref{sec:weak-form},
$S$ is expanded in the basis $\{\psi_j\}$, {so only the residual components in the range of $M_S$ \eqref{eq:Ms} are removed by the projector $\mathsf P_S$ in \eqref{eq:discrete-static-projector}.}
Therefore,
the discrete analogue of
Th.~\ref{thm:profiling}  holds only if this restricted static
forcing can absorb the difference between the unknown- and known-$\Ast$
systems, i.e., if the columns of $M_\eta-M_a$ lie in the range of $M_S$.  The
next lemma shows that this is the case for the spatial self-testing setting.

\begin{lemma}
\label{lem:self-test-compatibility}
Assume the spatial self-testing setting of Sec.~\ref{sec:weak-form},
the same spatial basis for $\eta$ and $S$, and a common quadrature and linear gradient
reconstruction in the known and unknown systems.
With $\mathsf R$ the spatial mass matrix in \eqref{eq:mass-stiffness}, the assembled matrices satisfy
\begin{equation}
 M_\eta-M_a=-M_S\mathsf R^{-1}\mathsf N, \quad \mathsf N_{pj}=\int_\Om\psi_j\nabla\Ast\cdot\nabla\psi_p\,\dd\bx.
 \label{eq:self-test-compatibility}
\end{equation}
Consequently,
$M_\eta^\top\mathsf Q^{-1}(I-\mathsf P_S)M_\eta\preceq M_a^\top\mathsf Q^{-1}M_a$. So, 
as in Th.~\ref{thm:profiling}, {information on $\eta$ may be lost when $\Ast$ is unknown}.
\end{lemma}

\begin{proof}
By \eqref{eq:MK} and \eqref{eq:Meta}, and since $\Ast$ does not depend on
$t$, the entry of $M_\eta-M_a$ in column $j$ and in the row of the test function
$\theta_m\psi_p$ of trajectory $i$ is
\[
 \frac1{NT}\Bigl(\int_0^T\theta_m\,\dd t\Bigr)
 \int_\Om\psi_j\nabla\Ast\cdot\nabla\psi_p\,\dd\bx
 =\frac1{NT}\Bigl(\int_0^T\theta_m\,\dd t\Bigr)\mathsf N_{pj}.
\]
By \eqref{eq:Ms}, the corresponding entry of $M_S$ is
$-\frac1{NT}\bigl(\int_0^T\theta_m\,\dd t\bigr)\mathsf R_{pj}$.  Hence, in each
group of rows with fixed $i$ and $m$, $M_\eta-M_a$ is obtained from $M_S$ by
replacing $-\mathsf R$ with $\mathsf N$, which is
\eqref{eq:self-test-compatibility}.  The same computation holds for the
assembled matrices when the same quadrature is used for all of them.

By \eqref{eq:self-test-compatibility}, every column of $M_\eta-M_a$ lies in
the range of $M_S$, which $I-\mathsf P_S$ maps to zero.  Therefore,
$(I-\mathsf P_S)M_\eta=(I-\mathsf P_S)M_a$: once the static forcing {has been eliminated,} the unknown-$\Ast$ system behaves exactly like the known-$\Ast$
system with the static-forcing directions removed.  Since $\mathsf P_S$ is a
$\mathsf Q^{-1}$-orthogonal projector, for every coefficient vector $\bc$,
\[
 \bc^\top M_\eta^\top\mathsf Q^{-1}(I-\mathsf P_S)M_\eta\bc
 =\norm{(I-\mathsf P_S)M_a\bc}_{\mathsf Q^{-1}}^2
 \leq\norm{M_a\bc}_{\mathsf Q^{-1}}^2
 =\bc^\top M_a^\top\mathsf Q^{-1}M_a\bc,
\]
where the inequality is the Pythagorean theorem. So, removing the component
$\mathsf P_SM_a\bc$ cannot increase the norm.
\end{proof}

\subsection{From a continuum rank condition to discrete solvability}

Th.~\ref{thm:identifiability} is a statement about the PDE problem, whereas
Prop.~\ref{prop:exact-discrete-solvability} is a statement about the
matrices of the learner.  This subsection connects the two: we ask whether
identifiability at the PDE level implies that the least-squares matrices have
full column rank.  The answer is yes if sufficiently many test functions are
used (Cor.~\ref{cor:discrete-solvability}) and for the spatial
self-testing design of Sec.~\ref{sec:weak-form}, provided its test
functions retain enough of the residual
(Prop.~\ref{prop:self-test-resolution}).

Testing a residual
$F\in Y$ against the test functions $\varphi_\ell$ of every trajectory gives
a vector $\br\in\R^{NK}$. This is how the right-hand sides and the columns of
the matrices in Sec.~\ref{sec:weak-form} are formed.  Let
$\Pi_{N,K}$ be the orthogonal projector in $Y$ onto the residuals that
correspond, through the Riesz map, to these test functions.  Then
\begin{equation}
 \norm{\Pi_{N,K}F}_Y^2=\br^\top\mathsf Q^{-1}\br,
 \label{eq:discrete-beta}
\end{equation}
that is, the discrete norm used by the learner is exactly the $Y$-norm of the
part of the residual seen by the test functions.

The first result shows that the rank condition of
Th.~\ref{thm:identifiability} yields full column rank once the test
space is rich enough. 

\begin{corollary}
\label{cor:discrete-solvability}
Assume the hypotheses of Th.~\ref{thm:identifiability} and, for each
trajectory $i$,
\begin{equation}
 t\longmapsto\nabla A_{(i)}(t,\cdot)
 \quad\text{belongs to }C([0,T];L^2(\Om)^2).
 \label{eq:time-continuity}
\end{equation}
Then, the map $\Tcal_U$ defined in \eqref{eq:TU-centered} is injective, so that for the diffusivity space $X_J$ spanned
by $\{\psi_j\}$,
\begin{equation}
 \inf_{0\neq u\in X_J}
 \frac{\norm{\Tcal_Uu}_Y}{\norm{u}_{L^2}}>0.
 \label{eq:continuous-discrete-beta}
\end{equation}
Moreover, suppose that the test spaces $V_K$ of Sec.~\ref{sec:weak-form}
are nested and that their union is dense in $L^2(0,T;H^1(\Om))$ as
$K\to\infty$, and keep the parameter spaces for
$\eta$, $S$, and $q$ fixed.  Then, for all sufficiently large $K$, the
matrices $M_a$, $[M_\eta\;M_S]$, and $M_q$ have full column rank, so the
conditions \eqref{eq:HK-discrete}, \eqref{eq:HU-discrete}, and
\eqref{eq:Hq-discrete} of Prop.~\ref{prop:exact-discrete-solvability}
hold.
\end{corollary}

\begin{proof}
If $\Tcal_Uu=0$ in $Y$, then $\Tcal_Ku$ is independent
of the trajectory and time. The time-continuity assumption extends this
equality to every selected trajectory--time
pair. For each selected quadruple, subtracting
the equation at $(i_0,t_0)$ from those at
$(i_k,t_k)$ gives
$\diver(u\nabla\delta A_k)=0$, $k=1,2,3$.
Lem.~\ref{lem:weak-cofactor}, together with
the assumed rank condition, then gives $u=0$.  This proves injectivity, and
the positive infimum in \eqref{eq:continuous-discrete-beta} follows by
compactness of the unit sphere in $X_J$.

The map $u\mapsto\Tcal_Ku$ is injective by the same three differences.
For the map $(u,s)\mapsto-\diver(u\nabla A_{(i)}(t,\cdot))-s$, taking differences eliminates $s$,
which does not depend on $i$ or $t$, and gives $u=0$. Then, the residual reduces
to $-s$, so $s=0$.  Finally, the source map is injective, since $q$ enters
every residual of the $\rho$ equation unchanged.  Density of the nested test
spaces gives $\Pi_{N,K}F\to F$ in $Y$, uniformly on the image of the unit
sphere of each fixed finite-dimensional parameter space.  Hence, each
projected map is injective for all sufficiently large $K$.  By
\eqref{eq:discrete-beta}, a projected map vanishes exactly when the
corresponding matrix--vector product vanishes, which proves full column rank.
\end{proof}

Cor.~\ref{cor:discrete-solvability} holds for each fixed parameter space as the test space is
enriched. However, it gives no bound that is uniform as the number of parameter basis
functions grows (see Sec.~\ref{sec:boundary-information}).  Moreover, the
spatial self-testing setting does not satisfy its density assumption:
increasing the number $M_t$ of temporal test functions does not enrich the
spatial tests, which remain tied to the parameter basis.
For the spatial self-testing setting, we therefore measure instead how much
of the residual its test functions retain. 

\begin{proposition}
\label{prop:self-test-resolution}
Consider the spatial self-testing setting, and let $\Pi_{N,K}$ be the
corresponding projector, with $K=M_tJ$.  Suppose that
\[
 \gamma:=\inf\frac{\norm{\Pi_{N,K}F}_Y}{\norm{F}_Y}>0,
\]
where the infimum is over all nonzero residuals with 
$F_{(i)}(t,\cdot)=-\diver\bigl(u\nabla A_{(i)}(t,\cdot)\bigr)-s$, $u,s\in X_J$.
Then, for every $u\in X_J$,
\begin{equation}
 \inf_{s\in X_J}\norm{\Pi_{N,K}\bigl(-\diver(u\nabla A_{(i)}(t,\cdot))-s\bigr)}_Y
 \geq\gamma\norm{\Tcal_Uu}_Y,
 \label{eq:self-test-infsup}
\end{equation}
where, as for $\Tcal_K$ and $\Tcal_U$, the residual inside the norm is
understood componentwise, for each trajectory $i$ and time $t$.
That is, when $\Ast$ is unknown, the information about $\eta$ retained by the
discrete learner is at least a fraction $\gamma^2$ of the PDE-level
information in Th.~\ref{thm:profiling}. Together with
\eqref{eq:continuous-discrete-beta}, the left-hand side is bounded below by a
positive multiple of $\norm{u}_{L^2}$.  The same holds when $\Ast$ is
known, with $\gamma$ defined over the residuals $\Tcal_Ku$ and $\Tcal_U$
replaced by $\Tcal_K$.
\end{proposition}

\begin{proof}
For every $u\in X_J$,
\begin{align*}
 \inf_{s\in X_J}\norm{\Pi_{N,K}\bigl(-\diver(u\nabla A_{(i)}(t,\cdot))-s\bigr)}_Y
 &\geq\gamma\inf_{s\in X_J}\norm{-\diver(u\nabla A_{(i)}(t,\cdot))-s}_Y\\
 &\geq\gamma\inf_{s\in H^{-1}(\Om)}\norm{-\diver(u\nabla A_{(i)}(t,\cdot))-s}_Y\\
 &=\gamma\norm{\Tcal_Uu}_Y.
\end{align*}
The first inequality is the definition of $\gamma$.  The second holds
because $X_J\subset H^{-1}(\Om)$, and minimizing over a larger set cannot
increase the minimum.  The equality is Th.~\ref{thm:profiling}:
$-\diver(u\nabla A_{(i)}(t,\cdot))$ differs from $\Tcal_Ku$ only by a term
independent of trajectory $i$ and time $t$, which the minimization over $s$ absorbs.
\end{proof}

\section{Boundary-induced loss of information}
\label{sec:boundary-information}

Th.~\ref{thm:identifiability} guarantees identifiability where the matrix
$\mathsf G$ has rank three.  
{However, the no-flux conditions force the normal derivatives of the observed fields to vanish at the boundary. As a result, the data carry less information near the boundary about the diffusivity $\eta$, which enters the equation through the gradient of the observed attractiveness.}
The source $q$ and the static forcing $S$ enter
the weak equations without spatial derivatives
of their coefficients. Thus, the normal-gradient
information-loss mechanism studied below is
specific to $\eta$ and does not
apply to $q$ or $S$. This section quantifies that loss of information on $\eta$: first locally
(in Sec.~\ref{sec:normal-info}), then for the rank condition of
Th.~\ref{thm:identifiability} (in Sec.~\ref{sec:rank-boundary}), and
finally for the smallest eigenvalue of the learner as the parameter space is
refined (in Sec.~\ref{sec:refinement-rate}).


\subsection{Normal information near a smooth boundary}
\label{sec:normal-info}

By Th.~\ref{thm:profiling}, a perturbation $u$ of $\eta$ changes the
residual through $u\nabla(A_{(i)}(t,\cdot)-\Ast)$ when $\Ast$ is known, and through
$u\nabla(A_{(i)}(t,\cdot)-\bar A)$ when $\Ast$ is unknown. 
Therefore, 
the information about $\eta$
available at a point $\bx$ in a unit direction $\boldsymbol\xi$ is measured by
the averaged squared component of these gradients in that direction. By introducing 
\begin{align}
 C_K(\bx)&:=\tavg{\nabla(A_{(i)}-\Ast)(t,\bx)\,\nabla(A_{(i)}-\Ast)(t,\bx)^\top},
 \label{eq:CK}\\
 C_U(\bx)&:=\tavg{\nabla(A_{(i)}-\bar A)(t,\bx)\,\nabla(A_{(i)}-\bar A)(t,\bx)^\top},
 \label{eq:CU}
\end{align}
we can write these measures as
$\boldsymbol\xi^\top C_K(\bx)\boldsymbol\xi$ and $\boldsymbol\xi^\top C_U(\bx)\boldsymbol\xi$ for known and unknown $\Ast$,
respectively. 

The next theorem shows that, near the boundary, the information on $\eta$ decays 
quadratically in the normal direction.


\begin{theorem}
\label{thm:boundary-loss}
Let $\bx$ be a point of a $C^2$ portion of $\partial\Om$, away from any
corner if $\Om$ is polygonal. Let $\mathbf n$ be the outward unit normal at
$\bx$, and let $\bx_\delta=\bx-\delta \mathbf n$ for small $\delta>0$.
Suppose that $A_{(i)}(t,\cdot)-\Ast$ is $C^2$ near $\bx$, uniformly in $i$
and $t$, and let $\kappa_{(i)}(t):=\partial^2_{\mathbf n\mathbf n}(A_{(i)}-\Ast)(t,\bx)$. 
Then, as $\delta\to0$,
\begin{align}
 \mathbf n^\top C_K(\bx_\delta)\mathbf n&=\delta^2\,\tavg{\kappa_{(i)}(t)^2}+o(\delta^2),
 \label{eq:known-boundary-law}\\
 \mathbf n^\top C_U(\bx_\delta)\mathbf n&=\delta^2\,\tavg{(\kappa_{(i)}(t)-\bar\kappa)^2}+o(\delta^2),
 ~ \bar\kappa:=\tavg{\kappa_{(i)}(t)}.
 \label{eq:unknown-boundary-law}
\end{align}
In particular, the information in the normal direction decays exactly like
$\delta^2$ when $\Ast$ is known, unless $\kappa_{(i)}(t)=0$ for every $i$ and almost
every $t$, and when $\Ast$ is unknown, unless $\kappa_{(i)}(t)$ is the same for every
$i$ and almost every $t$.
\end{theorem}

\begin{proof}
The no-flux conditions make the normal derivative of $A_{(i)}(t,\cdot)-\Ast$ vanish at
$\bx$. A Taylor expansion along the normal therefore gives
\begin{equation}
 \partial_{\mathbf n}(A_{(i)}-\Ast)(t,\bx_\delta)
 =-\delta\kappa_{(i)}(t)+o(\delta),
 \label{eq:normal-gradient-expansion}
\end{equation}
uniformly in $i$ and $t$.  Inserting \eqref{eq:normal-gradient-expansion} into
\eqref{eq:CK} gives
\[
 \mathbf n^\top C_K(\bx_\delta)\mathbf n
 = \frac1{NT} \left[\delta^2\tint{\kappa_{(i)}(t)^2}
 -2\delta\tint{\kappa_{(i)}(t)\,o(\delta)}
 +\tint{o(\delta)^2} \right],
\]
and the last two terms are $o(\delta^2)$. In fact, by the uniform $C^2$
assumption, the second derivatives of $A_{(i)}(t,\cdot)-\Ast$, and hence $\kappa_{(i)}(t)$, are
bounded independently of $i$ and $t$.  For $C_U$,
the normal derivative of $A_{(i)}(t,\cdot)-\bar A$ is that of $A_{(i)}(t,\cdot)-\Ast$ minus its average
over trajectories and times, so $\kappa_{(i)}(t)$ is replaced by
$\kappa_{(i)}(t)-\bar\kappa$.
\end{proof}

The uniform $C^2$ assumption can be relaxed to
requiring \eqref{eq:normal-gradient-expansion}
to hold in $L^2(0,T)^N$ for some
$\kappa\in L^2(0,T)^N$.


\subsection{The rank condition near the boundary}
\label{sec:rank-boundary}

Th.~\ref{thm:identifiability} requires the matrix $\mathsf G$ in
\eqref{eq:rank3matrix} to have rank three almost everywhere.  Its first two
columns contain the components of the gradients of the increments
$\delta A_k$.  At a boundary point, the no-flux conditions make the normal
components of these gradients vanish, so $\mathsf G$ loses rank there.  The
next corollary makes this precise: $\mathsf G$ has rank at most two on the
boundary and at distance $\delta$ from it, its smallest singular value grows
only linearly in $\delta$.  Since the smallest singular value measures how far
$\mathsf G$ is from losing rank, the rank condition of
Th.~\ref{thm:identifiability}, even when it holds almost everywhere,
becomes weaker and weaker as the boundary is approached.

\begin{corollary}
\label{cor:rank-boundary}
Let $\bx$, $\mathbf n$, and $\bx_\delta$ be as in
Th.~\ref{thm:boundary-loss}. Let $\bm\tau$ be a unit tangent at $\bx$,
and suppose the increments $\delta A_k$ in \eqref{eq:increments} are
$C^{2,1}$ near $\bx$.  Then, $\mathsf G(\bx)$ in 
\eqref{eq:rank3matrix}
has rank at most two and
$\sigma_{\min}(\mathsf G(\bx_\delta))=O(\delta)$ as $\delta\to0$.  In
addition, if the vectors
\begin{equation}
 \mathbf a=(\partial_{\bm\tau}\delta A_k)_{k=1}^3,
 \qquad
 \mathbf c=(\partial^2_{\mathbf n\mathbf n}\delta A_k)_{k=1}^3,
 \qquad
 \mathbf d=(\partial^2_{\bm\tau \bm\tau}\delta A_k)_{k=1}^3 ,
 \label{eq:acd-vectors}
\end{equation}
evaluated at $\bx$ are linearly independent, then
$\sigma_{\min}(\mathsf G(\bx_\delta))$ is bounded below by a positive
multiple of $\delta$, i.e., it decays exactly linearly.
\end{corollary}

\begin{proof}
By the no-flux conditions, $\partial_{\mathbf n}\delta A_k=0$ at $\bx$ for
$k=1,2,3$.  Since the first two columns of $\mathsf G$ contain the gradients of
the increments, $\mathsf G(\bx)(n_1,n_2,0)^\top=0$, so $\mathsf G(\bx)$ has
rank at most two.

Rotating the first two columns of $\mathsf G$ into the frame
$(\bm\tau,\mathbf n)$ does not change its singular values, and the Laplacian
is invariant under this rotation, $\Delta=\partial^2_{\bm\tau \bm\tau}
+\partial^2_{\mathbf n\mathbf n}$.  A Taylor expansion from $\bx$ to
$\bx_\delta$ therefore gives, column by column,
\begin{equation}
 \mathsf G(\bx_\delta)\ \simeq\ 
 \bigl[\,\mathbf a+O(\delta)\ \ \ -\delta\mathbf c+O(\delta^2)\ \ \ 
 \mathbf c+\mathbf d+O(\delta)\,\bigr],
 \label{eq:G-boundary-expansion}
\end{equation}
where $\simeq$ denotes equality up to the rotation.  The second column is
$O(\delta)$, hence $\sigma_{\min}(\mathsf G(\bx_\delta))=O(\delta)$.

If $\mathbf a$, $\mathbf c$, $\mathbf d$ are linearly independent, then,
since the determinant is linear in each column,
\eqref{eq:G-boundary-expansion} gives
$|\det\mathsf G(\bx_\delta)|=\delta\,|\det[\mathbf a,\mathbf c,\mathbf d]|
+O(\delta^2)$.  The determinant is the product of the three singular values,
and the two largest are bounded because the entries of $\mathsf G$ are bounded.
Hence, $\sigma_{\min}(\mathsf G(\bx_\delta))$ is bounded below by a positive
multiple of $\delta$ for small $\delta$.
\end{proof}

At a polygonal vertex, especially a re-entrant corner, the outward normal is
 not uniquely defined and solutions may have corner singularities. 

\subsection{Effect of the boundary on the learner as the parameter mesh is refined}
\label{sec:refinement-rate}
Secs.~\ref{sec:normal-info} and~\ref{sec:rank-boundary} 
describe the data near a single boundary point.
In this subsection, we ask what these local results imply
for the learner.  Its stability is governed by the smallest eigenvalue, relative
to the mass matrix $\mathsf R$ in \eqref{eq:mass-stiffness}, of the matrices
$M_a^\top\mathsf Q^{-1}M_a$ and $M_\eta^\top\mathsf Q^{-1}(I-\mathsf P_S)M_\eta$
in \eqref{eq:HK-discrete}--\eqref{eq:HU-discrete}.  
As the parameter mesh is refined, basis functions can
concentrate near the boundary, where the information is small by
Th.~\ref{thm:boundary-loss}, so this eigenvalue is expected to decrease.
Th.~\ref{thm:boundary-discrete-rate} below shows that it decreases like $1/J$,
where $J$ is the number of parameter basis functions, provided that information is lost only near the boundary.
In this subsection, we take the parameter space $X_J$ of
Sec.~\ref{sec:weak-form} to be a finite-element or spline space on a mesh of
size $h$, so that $h\sim J^{-1/2}$ in two dimensions.
Before stating Th.~\ref{thm:boundary-discrete-rate},
we need a geometric fact: on the
discrete parameter space, a function cannot be concentrated arbitrarily close
to the boundary.

\begin{lemma}
\label{lem:discrete-boundary-weight}
Let $\Om$ have $C^2$ boundary, and suppose that the mesh of $X_J$ is
boundary-fitted, shape-regular, and quasi-uniform, and that the polynomial
degree of $X_J$ is fixed as $J\to\infty$.  Then, there is a constant $c>0$, independent
of $J$, such that
\begin{equation}
 \norm{\operatorname{dist}(\cdot,\partial\Om)\,u_J}_{L^2(\Om)}
 \geq c\,J^{-1/2}\norm{u_J}_{L^2(\Om)}
 \qquad\text{for }u_J\in X_J.
 \label{eq:discrete-weighted-hardy}
\end{equation}
\end{lemma}

\begin{proof}
On each mesh cell $E$, $u_J$ is a polynomial of fixed degree.  Map $E$ to a
fixed reference cell.  There, all norms of such polynomials are equivalent,
so a polynomial cannot have most of its mass arbitrarily close to one side of
the cell.  If $E$ touches the boundary, the distance to $\partial\Om$ on $E$
scales like $h$ times the distance to the corresponding side of the reference
cell, which gives
\[
 \int_E\operatorname{dist}(\bx,\partial\Om)^2|u_J|^2\,\dd \bx
 \geq c_0^2h^2\int_E|u_J|^2\,\dd \bx.
\]
If $E$ does not touch the boundary, the same estimate follows from
$\operatorname{dist}(\bx,\partial\Om)\geq c_0h$ on $E$.  Shape regularity makes
$c_0$ independent of $E$.  Summing over the cells and using $h\sim J^{-1/2}$,
which holds for a quasi-uniform mesh in two dimensions, proves
\eqref{eq:discrete-weighted-hardy}.
\end{proof}

Let $\Pi_{N,K}$ be the
projector of Prop.~\ref{prop:self-test-resolution}, and define the information
about a perturbation $u_J\in X_J$ of $\eta$ in the settings of known and unknown $\Ast$ by
\begin{align*}
 \mathcal I_K(u_J)
 &:=\norm{\Pi_{N,K}\Tcal_Ku_J}_Y^2,\\
 \mathcal I_U(u_J)
 &:=\inf_{s_J\in X_J}
 \norm{\Pi_{N,K}\bigl(-\diver(u_J\nabla A_{(i)}(t,\cdot))-s_J\bigr)}_Y^2.
\end{align*}
In terms of coefficient vectors, these are the quadratic forms of the
matrices $M_a^\top\mathsf Q^{-1}M_a$ and
$M_\eta^\top\mathsf Q^{-1}(I-\mathsf P_S)M_\eta$ of
Prop.~\ref{prop:exact-discrete-solvability}.  In the next theorem, $\mathcal I$ denotes either
$\mathcal I_K$ or $\mathcal I_U$, with the assumptions imposed on that choice.  Finally, set
\begin{equation}
 \lambda_{\min,J}
 :=\inf_{0\neq u_J\in X_J}
 \frac{\mathcal I(u_J)}{\norm{u_J}_{L^2}^2}.
 \label{eq:discrete-information-rate-definitions}
\end{equation}

\begin{theorem}
\label{thm:boundary-discrete-rate}
Assume the setting of Lem.~\ref{lem:discrete-boundary-weight}.  Suppose the
parameter and test spaces are refined compatibly and the resolved information
form has the following two uniform properties:
\begin{align}
 \mathcal I(u_J)
 &\geq c_I\norm{\operatorname{dist}(\cdot,\partial\Om)\,u_J}_{L^2}^2
 ~\text{for every }u_J\in X_J,
 \label{eq:weighted-information-coercivity}\\
 \mathcal I(w_J)
 &\leq C_I\norm{\operatorname{dist}(\cdot,\partial\Om)\,w_J}_{L^2}^2,
 ~
 \norm{\operatorname{dist}(\cdot,\partial\Om)\,w_J}_{L^2}
 \leq C_BJ^{-1/2}\norm{w_J}_{L^2}
 \label{eq:normal-boundary-trial}
\end{align}
for some $0\neq w_J\in X_J$.
Here, $w_J$ is a parameter perturbation supported near a smooth boundary
segment that varies mainly in the normal direction, i.e., in the direction in
which, by Th.~\ref{thm:boundary-loss}, the observed fields carry little
information.  Then, with $c$ the constant of
Lem.~\ref{lem:discrete-boundary-weight},
\begin{equation}
 {
 \frac{c_Ic^2}{J}\leq\lambda_{\min,J}\leq\frac{C_IC_B^2}{J}.}
 \label{eq:two-sided-J-rate}
\end{equation}
\end{theorem}

\begin{proof}
Combining \eqref{eq:weighted-information-coercivity} with
Lem.~\ref{lem:discrete-boundary-weight} gives, for every $u_J\in X_J$,
\[
 \mathcal I(u_J)
 \geq\frac{c_Ic^2}{J}\norm{u_J}_{L^2}^2.
\]
Taking the infimum gives the lower bound for $\lambda_{\min,J}$.  Testing its
Rayleigh quotient with the boundary-layer function in
\eqref{eq:normal-boundary-trial} gives
\[
 \lambda_{\min,J}
 \leq\frac{\mathcal I(w_J)}{\norm{w_J}_{L^2}^2}
 \leq\frac{C_IC_B^2}{J}.
\]
\end{proof}

Th.~\ref{thm:boundary-discrete-rate} says that, as the parameter mesh is
refined, the learner loses stability at a definite rate: the smallest
eigenvalue $\lambda_{\min,J}$ decreases like $1/J$.
The upper bound reflects the boundary-induced loss
of information, while the lower bound ensures that
this loss is no worse than $1/J$.
The latter relies on the global coercivity condition
\eqref{eq:weighted-information-coercivity}, which
does not follow from the local boundary analysis in
Th.~\ref{thm:boundary-loss}. The numerical experiments in Sec.~\ref{sec:num_res}
examine this loss of stability under refinement.

\section{Numerical results}\label{sec:num_res}

We generate observations by finite element simulations of
\eqref{eq:PDE1}--\eqref{eq:PDE2}, implemented in FEniCSx \cite{OurCode}
and following the numerical framework
of \cite{Hao2026}. We use continuous $P_1$ finite elements
on triangular meshes, backward Euler time discretization,
and the iterative partitioned solver described therein. The time step is $\Delta t_{\rm sim}=0.005$,
and simulations run until $T_{\rm sim}=15$. We consider two configurations that we call Neumann square and Chicago. The geometry, mesh, and ground-truth
fields for each configuration are described below.

The Neumann square is the domain $\Omega=[0,L]^2$, with $L=2.8$.
We use a uniform triangular mesh obtained by dividing each
cell of a Cartesian grid of spacing $h_{\rm sim}=0.0125$
into two triangles, and we impose no-flux conditions
on all four sides.
The prescribed fields $\eta$, $q$, and $\Ast$ are smooth, with nominal ranges
$[0.02,0.05]$ for $\eta$, $[1.15,1.85]$ for $q$,
and $[1/30,0.19]$ for $\Ast$.
The static field $\Ast$ satisfies
\eqref{eq:nofluxAst} on every side of the square.
The explicit field definitions are given in App.~\ref{sec:app}.

For the Chicago experiment, we use a polygonal approximation
of the city boundary \cite{Chicago}, rescaled so that the longest side
of its bounding box is 5. See details in App.~\ref{sec:app}. 
An unstructured triangular mesh is generated using Gmsh \cite{geuzaine2009gmsh}
with target mesh size $h_{\rm sim}=0.0125$.
No-flux boundary conditions are imposed on the entire boundary. The diffusivity $\eta$ is prescribed as a rescaled combination
of Gaussian functions, while $q$ combines a linear spatial
trend with a Gaussian component.
Their nominal ranges are $[0.02,0.05]$ and $[1.15,1.85]$,
respectively.
The static field $\Ast$ is obtained by solving $\Leta\Ast=S$ (see \eqref{eq:static-source}) with boundary condition \eqref{eq:nofluxAst}, in the $P_1$ finite element space used for the simulations, where $S$ is a smooth Gaussian mixture with nominal range $[1/30,0.20]$.
The explicit definitions of $\eta$, $q$, and $S$ are given in App.~\ref{sec:app}. 


For each independent trajectory, the initial conditions are
\[
A_0=(q+\Ast)e^{\epsilon_0 g_A-\epsilon_0^2/2},
\qquad
\rho_0=\frac{q}{q+\Ast}
e^{\epsilon_0 g_\rho-\epsilon_0^2/2},
\qquad \epsilon_0=0.3,
\]
where $g_A$ and $g_\rho$ are independent, unit-pointwise-variance
Gaussian random fields constructed from band-limited Fourier
series, with correlation-length parameter $0.4$ and maximum
wavenumber $7.5$.

For both configurations, 
the simulated scalar fields are sampled on an observation grid
of spacing $h_{\mathrm{obs}}=0.025$ at time intervals
$\Delta t_{\rm obs}=0.01$. {The observations are sampled uniformly in space and time.
Nonuniform observations could also be used, with suitable
interpolation or numerical integration procedures.}
No noise is added to the observations, so that the experiments isolate the
effects of identifiability and discretization. The effect of observation
noise will be studied in Part~II.
Unless otherwise stated, learning uses $N = 4$ trajectories over
$[0,T_{\rm obs}]$ with $T_{\rm obs} = 8$.
The $P_1$ finite element space is distinct from the cubic
B-spline spaces used for parameter recovery.  

We recall that we assume to have direct observations of $A$ and $\rho$ only. Spatial gradients are reconstructed from the scalar observations
using moving least squares \cite{lancaster1981surfaces}.
At each observation point, we fit a local polynomial of total
degree three using $36$  neighboring points and evaluate its
derivatives at the center. 
{For the flux term in \eqref{eq:PDE2}, we use the identity
\begin{equation}
 \nabla\rho_{(i)}-2\rho_{(i)}\nabla\log A_{(i)}
 =A_{(i)}^2\nabla\!\left(\frac{\rho_{(i)}}{A_{(i)}^2}\right)
 \label{eq:structure-flux}
\end{equation}
to reconstruct the composite flux directly from the
observed scalar fields, rather than subtracting two
independently approximated gradients.}

In all the tests, we use $\mathsf Q^{-1}$-weighted least squares, without parameter regularization, 
as explained in Sec.~\ref{sec:weak-form}. We recall that we do not enforce positivity constraints, as specified in Rem.~\ref{rem:pos_eta}.


\subsection{Parameter recovery}\label{sec:param_rec}

Let us start by considering the Neumann square test. 
Each unknown field is represented by $9\times9$ cubic
tensor-product B-splines ($J=81$).
We use the spatial self-testing construction of
Sec.~\ref{sec:weak-form}, with
$\theta_m(t)=\sin({m\pi t}/{T_{\rm obs}})$, $m=1,\ldots,10$. 

\begin{figure}[htb!]
    \centering
    \captionsetup[subfigure]{position=top}
    \begin{subfigure}{0.49\linewidth}
        \caption{$q$ -- both known and unknown $\Ast$ cases}
        \includegraphics[width=\linewidth]{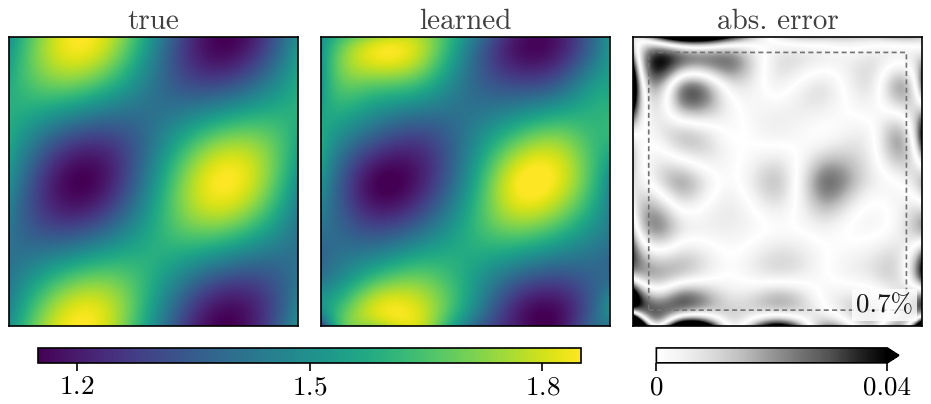}
        \label{fig:square-q}
    \end{subfigure}\hfill
    \begin{subfigure}{0.49\linewidth}
        \caption{$\Ast$ -- unknown $\Ast$ case}
        \includegraphics[width=\linewidth]{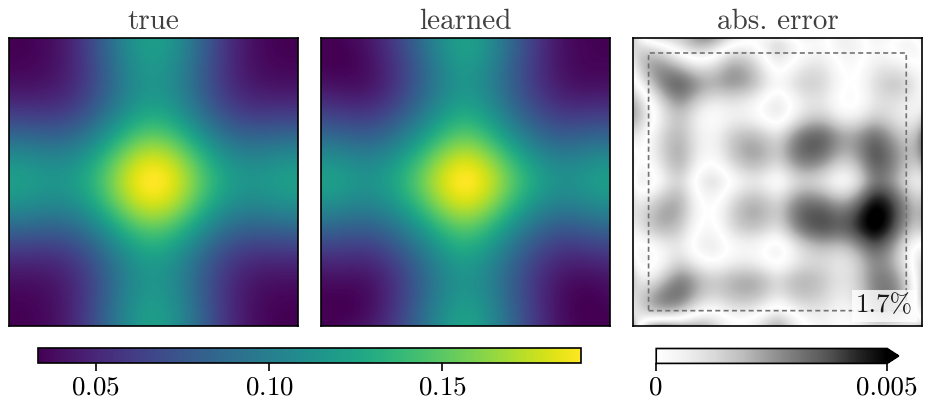}
        \label{fig:square-Ast-unknown}
    \end{subfigure}

    \vspace{1ex}
    \begin{subfigure}{0.49\linewidth}
        \caption{$\eta$ -- known $\Ast$ case}
        \includegraphics[width=\linewidth]{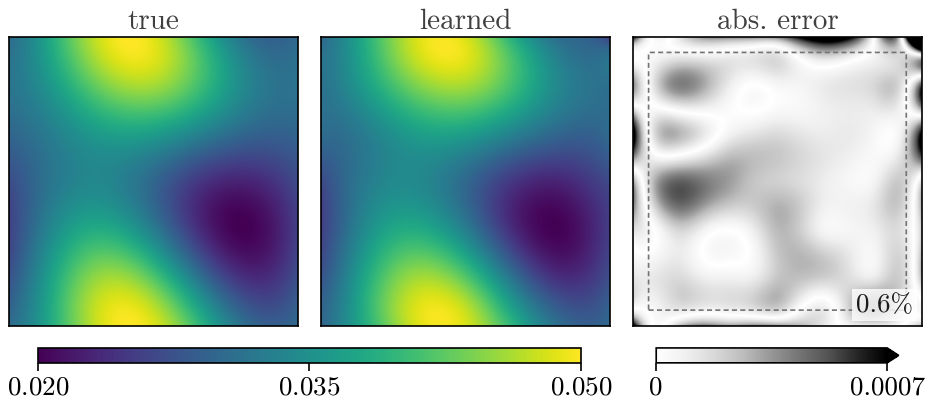}
        \label{fig:square-eta-known}
    \end{subfigure}\hfill
    \begin{subfigure}{0.49\linewidth}
        \caption{$\eta$ -- unknown $\Ast$ case}
        \includegraphics[width=\linewidth]{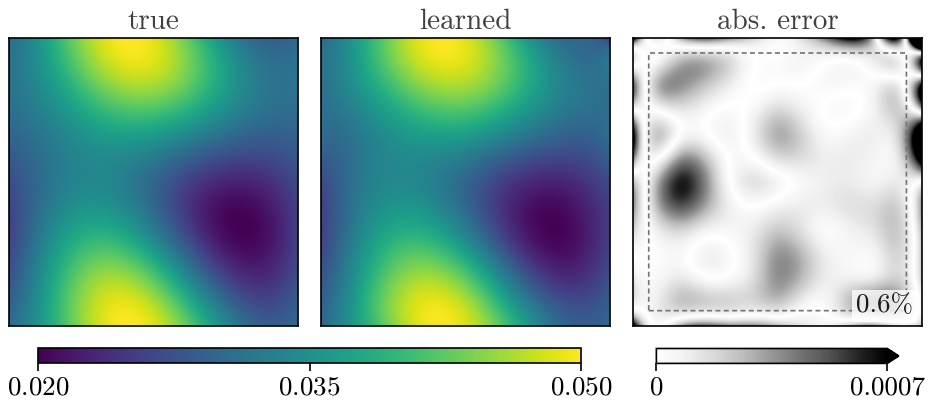}
        \label{fig:square-eta-unknown}
    \end{subfigure}

    \caption{Neumann square: (a) learning of 
    $q$ in the cases of both $\Ast$ known and unknown, 
    (b) learning of $\Ast$ when $\Ast$
    is unknown, and learning of 
    $\eta$ in the case of (c) $\Ast$ known and (d) $\Ast$ unknown. 
    Each panel displays the ground truth, reconstruction,
    and pointwise absolute error.
    Percentages in the error panels indicate
    full-domain relative $L^2$ errors. 
    The dashed line is at distance $\delta = 0.15$ from the boundary. 
    }
    \label{fig:square-recovery}
\end{figure}

Fig.~\ref{fig:square-recovery} shows
the recovery of $q$ and $\eta$ for both known and unknown $\Ast$ and the recovery of $\Ast$ when $\Ast$ is unknown.
Our weak-form learner
recovers all three spatially varying parameters accurately.  In every case, the learned field is visually
indistinguishable from the ground truth: the locations, shapes, and
amplitudes of the maxima and minima are all reproduced.  These results
are obtained by using the four trajectories together: their residuals
enter a single least-squares problem, \eqref{eq:weighted-learner_AstK}
or \eqref{eq:weighted-known-q-learner_AstU} with $N=4$.  We refer to
this as the {joint fit}.  The joint fit
yields relative $L^2$ errors of approximately $0.7\%$ for $q$
(Fig.~\ref{fig:square-q}), $0.6\%$ for $\eta$
(Figs.~\ref{fig:square-eta-known} and~\ref{fig:square-eta-unknown}),
and $1.7\%$ for $\Ast$ (Fig.~\ref{fig:square-Ast-unknown}).

The pointwise errors in Fig.~\ref{fig:square-recovery} are small compared to the ranges of the
corresponding fields.  
For both $q$ and $\eta$, the largest errors
are located near the boundary of the domain.
For $\eta$, this is consistent with the
boundary-induced loss of information analyzed
in Sec.~\ref{sec:boundary-information}.
The boundary errors in $q$ are not explained
by the same mechanism and may instead be
related to discretization errors, including
the reconstruction of the flux in
\eqref{eq:structure-flux} near the boundary.
{Treating $\Ast$ as unknown does not noticeably affect the recovery of $\eta$ in the joint fit: the relative error remains at $0.6\%$ in both settings, and the error fields in Figs.~\ref{fig:square-eta-known} and~\ref{fig:square-eta-unknown} have comparable magnitudes, although their spatial patterns differ. This does not contradict Th.~\ref{thm:profiling}, which concerns the information about $\eta$, not the error. The information is indeed smaller when $\Ast$ is unknown (see Sec.~\ref{sec:info_loss}), but with noise-free data and four trajectories, this loss has little effect on the error, which might be dominated by the approximation of the parameters and the reconstruction of the gradients.}  
The static attractiveness has the largest error, $1.7\%$, with
absolute errors up to about $0.005$ for values in $[1/30,0.19]$.
This is expected, since $\Ast$ is not learned directly: it is
recovered from the learned $\eta_J$ and $S_J$ through the elliptic
problem \eqref{eq:recover-Ast}, so the errors in both propagate to
$\Ast$.  Even so, the reconstruction captures the central maximum and
the minima at the four corners of the domain.

To assess the benefit of the joint fit, we repeat the learning with
each of the four trajectories separately, while keeping all other settings
fixed.  Tab.~\ref{tab:single-vs-joint} compares the error of the
joint fit with the mean and standard deviation of the errors of the
four single-trajectory fits.  The joint fit is more accurate for every
parameter: it reduces the error by a factor of about two for $q$ and
$\Ast$ and a factor of about three (resp., five) for $\eta$ when $\Ast$ is known (resp., unknown). The largest gain occurs
where the most information is lost, i.e., when $\eta$ is learned
without knowledge of $\Ast$.  With a single trajectory, not knowing
$\Ast$ increases both the mean error for $\eta$  and its variability across trajectories, so the result depends heavily on which
trajectory is observed.  
This agrees with the
analysis of Sec.~\ref{sec:ident}: when $\Ast$ is unknown, only the
variation of the observed attractiveness about its average carries
information about $\eta$ (see \eqref{eq:TU-centered}), and
independent trajectories add precisely this kind of variation.

\begin{table}[htb!]
\centering
\small
\setlength{\tabcolsep}{4pt}
\begin{tabular}{ccc}
\toprule
Parameter 
& \shortstack{Joint fit}
& \shortstack{Single traj. fit (mean $\pm$ std)} \\
\midrule
$q$ - known and unknown $\Ast$
& 0.7 & $1.5 \pm 0.3$ \\
$\eta$ - known $\Ast$
 & 0.6 & $2.1 \pm 0.5$ \\
$\eta$ - unknown $\Ast$
& 0.6 & $3.2 \pm 1.5$ \\
$\Ast$ - unknown $\Ast$
 & 1.7 & $3.6 \pm 1.4$ \\
\bottomrule
\end{tabular}
    \caption{
    Neumann square: 
    full-domain relative $L^2$ errors (\%) in the case of 
    joint fit and single-trajectory fit. 
    The single-trajectory column reports the mean and sample
    standard deviation of four separate fits, one per trajectory.}
    \label{tab:single-vs-joint}
\end{table}

\begin{figure}[htb!]
    \centering
    \includegraphics[width=\linewidth]{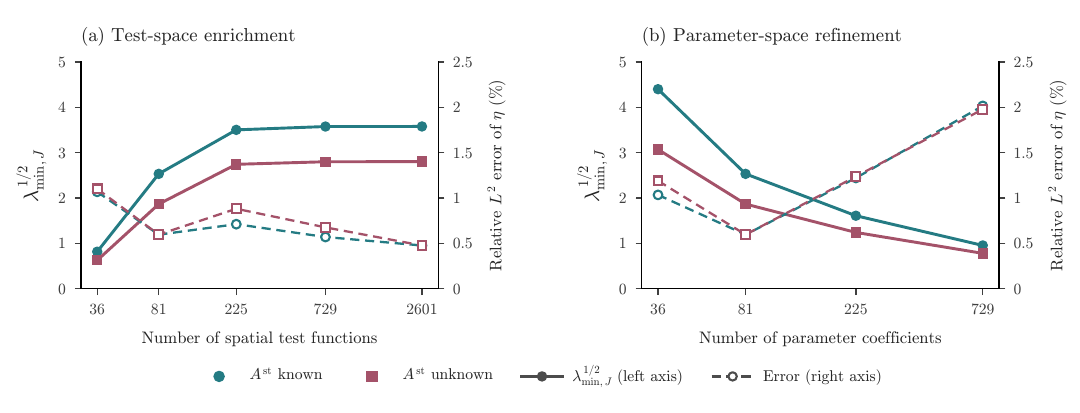}
    \caption{Neumann square: $\lambda_{\min,J}^{1/2}$
\eqref{eq:discrete-information-rate-definitions} in solid line and relative $L^2$ errors of $\eta$ in dashed line as  (a) the number of spatial test functions increases when using 81 parameter coefficients (81 spatial tests correspond to self-testing)
and 
(b) the number of parameter coefficients increases under self-testing. The left vertical axis corresponds to $\lambda_{\min,J}^{1/2}$, while the right vertical axis refers to the relative $L^2$ errors of $\eta$.}
    \label{fig:spaces}
\end{figure}

Next, we justify the choice of $J=81$ coefficients and spatial
self-testing.  As a measure of the stability of the learner, we use $\lambda_{\min,J}^{1/2}$, with $\lambda_{\min,J}$ the smallest eigenvalue of the learner defined in \eqref{eq:discrete-information-rate-definitions}.
In Fig.~\ref{fig:spaces}(a), we keep $J=81$ and enrich the spatial
test space beyond self-testing, which corresponds to 81 spatial test
functions.  The value of $\lambda_{\min,J}^{1/2}$ increases up to 225 test
functions and then levels off. 
However, the relative error of $\eta$
changes little beyond self-testing: it drops from about $1.1\%$ with
36 test functions to $0.6\%$ with 81, and remains between $0.5\%$ and
$0.9\%$ for richer test spaces.  Hence, self-testing already gives
nearly the smallest error, with far fewer test functions.

In Fig.~\ref{fig:spaces}(b), we refine the parameter and test spaces
together, keeping the self-testing setting.  As $J$ increases, $\lambda_{\min,J}^{1/2}$
decreases in both static-field settings.  This is the loss of
stability under parameter refinement described by
Th.~\ref{thm:boundary-discrete-rate}. 
The error of $\eta$ is smallest for $J=81$, at about $0.6\%$
in both settings.  It is larger for the coarsest space, $J=36$, and it
grows to about $2\%$ for $J=729$, while $\lambda_{\min,J}^{1/2}$ keeps
decreasing.  Among the tested resolutions, $J=81$ therefore gives the
smallest errors, which motivates our choice.

\begin{figure}[htb!]
    \centering
    \captionsetup[subfigure]{position=top}
    \begin{subfigure}{0.49\linewidth}
        \caption{$q$ -- both known and unknown $\Ast$}
        \includegraphics[width=\linewidth]{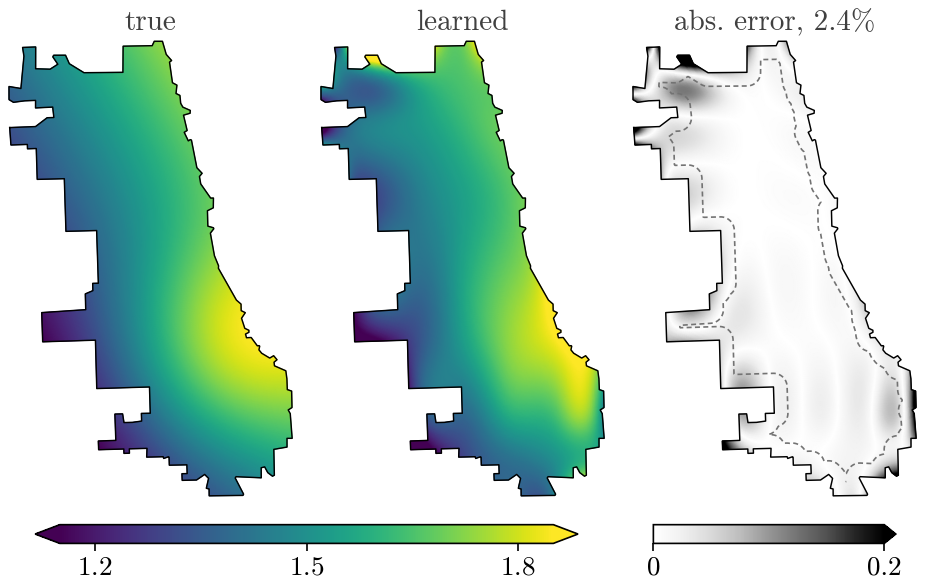}
        \label{fig:chicago-q}
    \end{subfigure}\hfill
    \begin{subfigure}{0.49\linewidth}
        \caption{$\Ast$ -- unknown $\Ast$}
        \includegraphics[width=\linewidth]{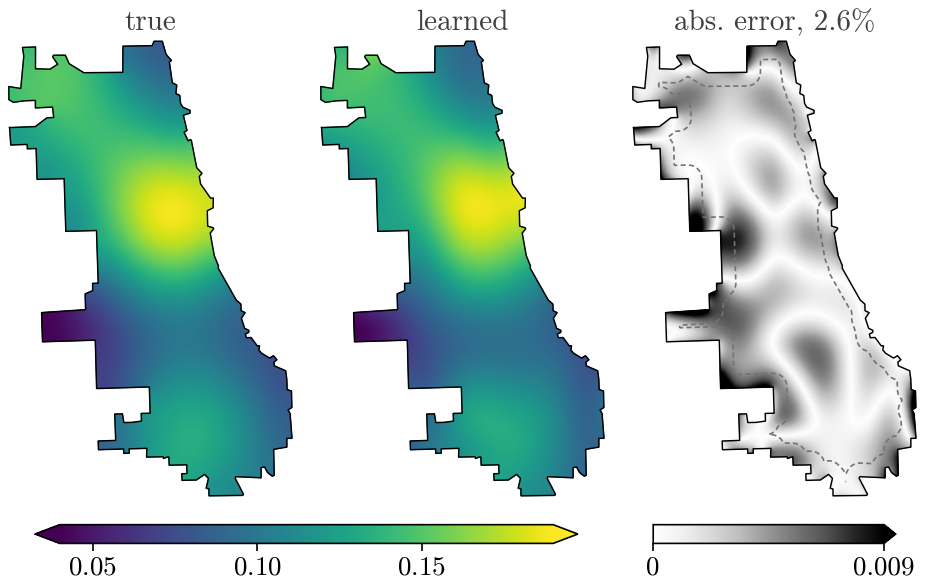}
        \label{fig:chicago-Ast-unknown}
    \end{subfigure}

    \vspace{1ex}
    \begin{subfigure}{0.49\linewidth}
        \caption{$\eta$ -- known $\Ast$}
        \includegraphics[width=\linewidth]{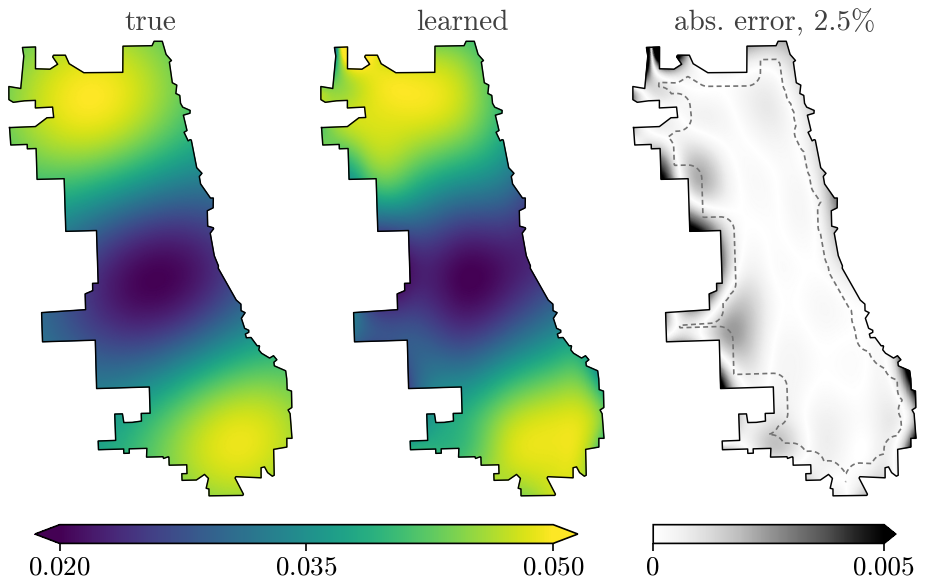}
        \label{fig:chicago-eta-known}
    \end{subfigure}\hfill
    \begin{subfigure}{0.49\linewidth}
        \caption{$\eta$ -- unknown $\Ast$}
        \includegraphics[width=\linewidth]{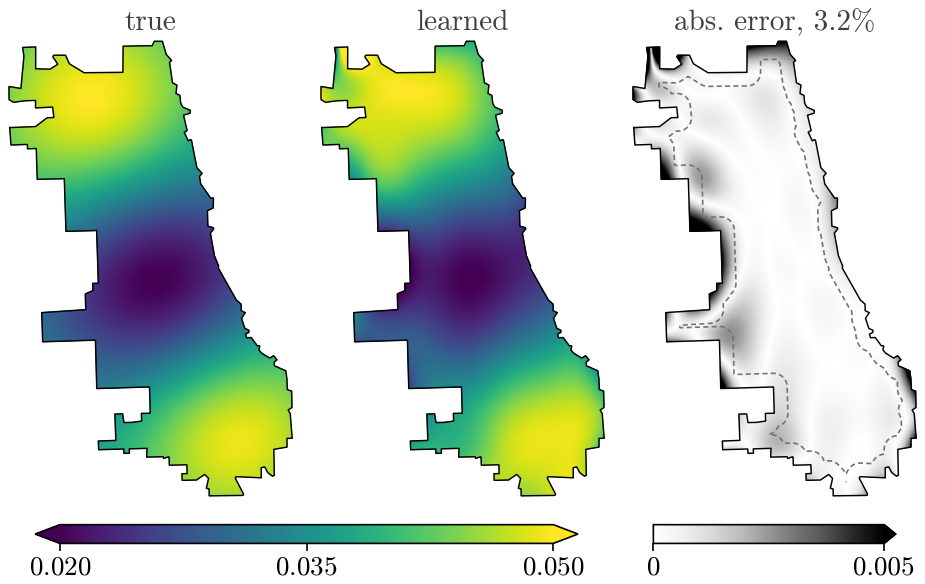}
        \label{fig:chicago-eta-unknown}
    \end{subfigure}

    \caption{Chicago: (a) learning of 
    $q$ in the cases of both $\Ast$ known and unknown, (b) learning of $\Ast$ when $\Ast$
    is unknown, and learning of 
    $\eta$ in the case of (c) $\Ast$ known and (d) $\Ast$ unknown. 
    Each panel displays the ground truth, reconstruction,
    and pointwise absolute error.
    Percentages in the error panels indicate
    full-domain relative $L^2$ errors. The dashed line is at distance $\delta = 0.15$ from the boundary.}
    \label{fig:chicago_recovery}
\end{figure}

For the Chicago test,
the cubic tensor-product B-splines are restricted to the domain,
with 59 active basis functions.
We retain spatial self-testing and the same temporal tests
and gradient-reconstruction settings as in the Neumann square experiment. As before, we use the joint fit of four trajectories. Fig.~\ref{fig:chicago_recovery} shows
the recovery of $q$ and $\eta$ for both known and unknown $\Ast$ and the recovery of $\Ast$ when $\Ast$ is unknown. We see that our method also recovers
the spatially varying parameters on this irregular domain.  The learned
fields reproduce all the main features of the ground truth: the
large-scale trend of $q$ and its maximum in the south-eastern part of
the city (Fig.~\ref{fig:chicago-q}), the two maxima of $\eta$ in the
north and south and its minimum in the center
(Figs.~\ref{fig:chicago-eta-known} and~\ref{fig:chicago-eta-unknown}),
and the location and shape of the maxima of $\Ast$
(Fig.~\ref{fig:chicago-Ast-unknown}).  The relative $L^2$ errors are
approximately $2.4\%$ for $q$, $2.5\%$ and $3.2\%$ for $\eta$ with
known and unknown $\Ast$, respectively, and $2.6\%$ for $\Ast$.  These
errors are larger than for the Neumann square, but they remain small for a challenging case.

The error panels in Fig.~\ref{fig:chicago_recovery}
show that, for both $q$ and $\eta$, the larger
errors are concentrated in a thin band along
the boundary, outside the dashed contour.
Note that the error color scales are capped,
so the largest errors along the boundary
exceed the maximum value shown. As for the Neumann square, the larger errors in $\eta$ are consistent with the boundary-induced loss of information analyzed in Sec.~\ref{sec:boundary-information}, which here applies to the edges of the polygonal boundary away from its corners, while those in $q$ may be related to the reconstruction of the composite flux \eqref{eq:structure-flux}, where the MLS neighborhoods become one-sided.
For this test,
unlike the Neumann square test, not knowing $\Ast$ has a visible cost: the error for $\eta$ increases from $2.5\%$ to $3.2\%$.  The
error remains localized near the boundary in both settings, so the
recovery of $\eta$ in the interior is essentially unaffected.  The
error for $\Ast$ is also largest near the boundary, but, unlike for $q$
and $\eta$, it is visible in the interior as well.  


\subsection{Information loss from an unknown $\Ast$ and effect of the observation window}\label{sec:info_loss}

In this subsection, we relate the data and the learner to the
analysis of Sec.~\ref{sec:ident}, and address two questions: how much
information about $\eta$ is lost when $\Ast$ is unknown, and how the
information depends on when the trajectories are observed.

We consider the Neumann square with the same four trajectories as in
Sec.~\ref{sec:param_rec}, but we restrict the observations to a window
$[t_0,t_0+2]$, which we slide along the trajectories.
The sampling
intervals, the parameter space, and the test functions are kept fixed
relative to the window.  For each window, we compute  two quantities.

\begin{figure}[htb!]
    \centering
    \includegraphics[width=\linewidth]{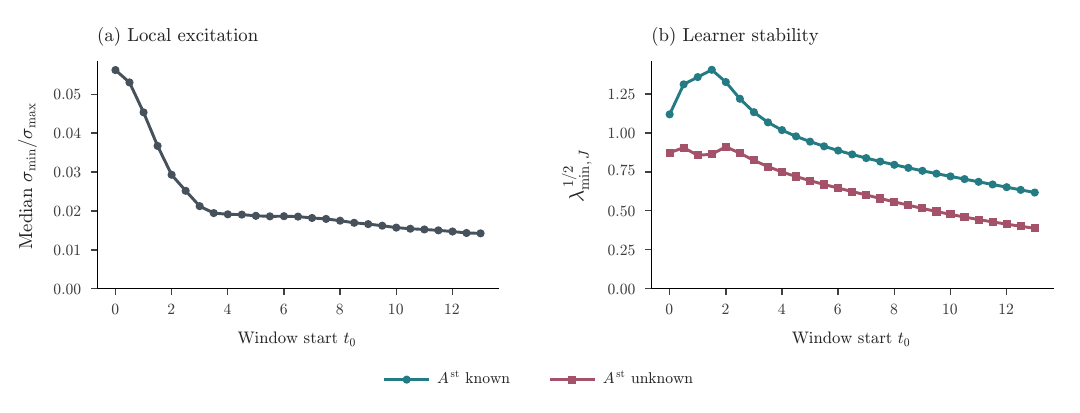}
    \caption{
    Neumann square:
    (a) spatial median over $\Om_{0.15}$ of $\sigma_{\min}(\widetilde{\mathsf G}(\bx))/\sigma_{\max}(\widetilde{\mathsf G}(\bx))$ and
    (b) $\lambda_{\min,J}^{1/2}$ for known and unknown $\Ast$.}
    \label{fig:window-identifiability}
\end{figure}

The first quantity concerns the data alone and extends the rank condition
\eqref{eq:rank3condition} to all observations in the window at once.  At each point $\bx$, we collect the vectors
$(\partial_{x_1}A_{(i)},\partial_{x_2}A_{(i)},\Delta A_{(i)})$ for all
trajectories and all sampled times in the window, subtract their weighted
mean (i.e., center them), and stack the resulting vectors as the rows of a matrix
$\widetilde{\mathsf G}(\bx)$.  The weights are those of the trapezoidal rule
in time, with equal weighting across trajectories, and each centered vector
is multiplied by the square root of its weight before being stacked.
Each window contains 201 sampled times per trajectory, so
$\widetilde{\mathsf G}(\bx)\in\R^{804\times3}$.  
We compute the ratio of its smallest to its largest singular value
$\sigma_{\min}(\widetilde{\mathsf G}(\bx))/\sigma_{\max}(\widetilde{\mathsf G}(\bx))$.
 $\widetilde{\mathsf G}(\bx)$ has rank three exactly when the observed attractiveness fields at some four trajectory--time pairs, possibly depending on $\bx$, give a rank-three increment matrix at $\bx$, which is the condition of Rem.~\ref{rem:multi-snapshot-rank}.  Subtracting the mean also reflects the setting of unknown $\Ast$, in which only the variation of the observations
around their average informs $\eta$ (Th.~\ref{thm:profiling}).  A ratio
close to zero means that the observed states are close to violating the
rank condition.  Since the rank condition degenerates near the boundary
(Cor.~\ref{cor:rank-boundary}), we report the spatial median of the
ratio over the interior region
$\Om_{0.15}=\{\bx\in\Om:\operatorname{dist}(\bx,\partial\Om)>0.15\}$.
Fig.~\ref{fig:window-identifiability}(a) shows how $\sigma_{\min}(\widetilde{\mathsf G}(\bx))/\sigma_{\max}(\widetilde{\mathsf G}(\bx))$ varies as $t_0$ is increased.
The second quantity is the stability measure $\lambda_{\min,J}^{1/2}$ of the
learner introduced in Sec.~\ref{sec:param_rec}, computed for known and
unknown $\Ast$. This is shown in Fig.~\ref{fig:window-identifiability}(b).

 The median
ratio in Fig.~\ref{fig:window-identifiability}(a) drops from about $0.056$
for $t_0=0$ to about $0.02$ for $t_0\approx3$, and then decreases slowly.
The value of $\lambda_{\min,J}^{1/2}$ in
Fig.~\ref{fig:window-identifiability}(b) is largest for $t_0$ between
$1.5$ and $2$, and then decreases steadily in both settings. This is
consistent with the decay of the transients: at later times the observed
fields vary less, and therefore constrain $\eta$ less.
Hence, observations
that capture the early transients are the most informative.  For the
earliest windows, the two quantities do not vary together: $\sigma_{\min}(\widetilde{\mathsf G}(\bx))/\sigma_{\max}(\widetilde{\mathsf G}(\bx))$
decreases while $\lambda_{\min,J}^{1/2}$ increases.  This is possible because
the ratio is a pointwise property of the data, whereas
$\lambda_{\min,J}$ measures the information retained by the discrete learner,
which also depends on the test functions and on the parameter space.
Fig.~\ref{fig:window-identifiability}(b) shows that
$\lambda_{\min,J}^{1/2}$ is smaller when $\Ast$ is unknown, for every window,
as predicted by Th.~\ref{thm:profiling} and
Lem.~\ref{lem:self-test-compatibility}.  The loss is substantial for every window.

\subsection{Boundary-induced information loss}\label{sec:bc-info-loss}

In this subsection, we test the predictions of Sec.~\ref{sec:boundary-information} on the loss of
information near a no-flux boundary. Locally, we
check whether the information contained in the data decays toward the
boundary at the rates given by Th.~\ref{thm:boundary-loss} and
Cor.~\ref{cor:rank-boundary}.  Globally, we check whether the
parameter directions that the learner constrains least are located near
the boundary, as suggested by Th.~\ref{thm:boundary-discrete-rate}.
In both cases, we compare the Neumann square with a periodic square, which
serves as a reference without a physical boundary.  The periodic square uses
the same domain, mesh, prescribed fields, and initial conditions as the
Neumann square, but periodic boundary conditions are imposed in both
coordinate directions instead of no-flux conditions.  This is possible
because the prescribed fields are periodic (see App.~\ref{sec:app}).  Since the periodic square has
no boundary, no boundary-induced loss of information is expected.
All experiments use four trajectories over $[0,8]$,
as in Sec.~\ref{sec:param_rec}.

\begin{figure}[htb!]
    \centering
    \includegraphics[width=\linewidth]{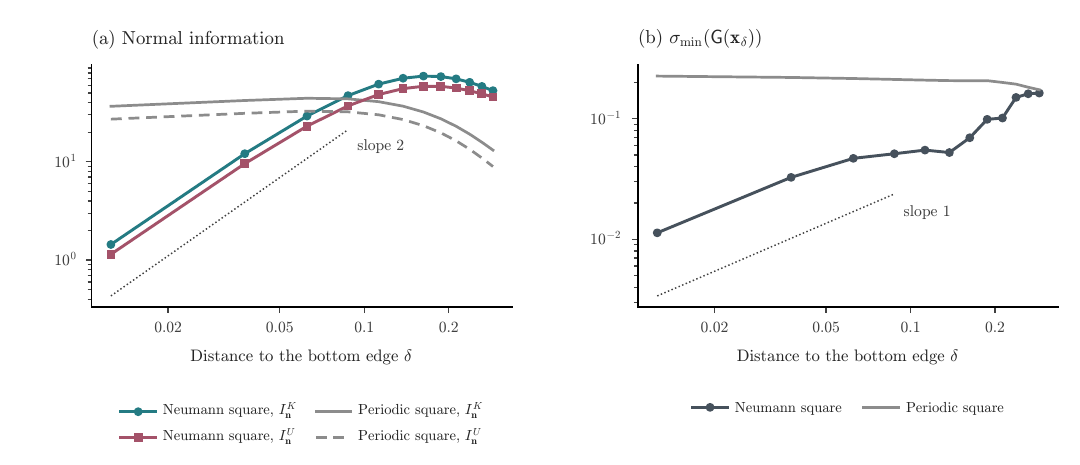}
    \caption{Neumann square: (a) normal information $I_{\mathbf n}^K$ and $I_{\mathbf n}^U$ for known
    and unknown $\Ast$ and (b) smallest singular value of the increment matrix
    $\mathsf G(\bx_\delta)$ as
    functions of the distance $\delta$ to the bottom edge.
    Grey curves show the same quantities for the periodic square, evaluated at the same coordinates.
    Dotted lines indicate reference slopes.}
    \label{fig:boundary-law}
\end{figure}

We sample the Neumann square along normals to the central half of its
bottom edge, $x_1\in[L/4,3L/4]$, and denote by $\delta=x_2$ the distance to
that edge.  Staying away from the corners ensures that the sampled points
lie near a straight portion of the boundary, where
Th.~\ref{thm:boundary-loss} and Cor.~\ref{cor:rank-boundary}
apply.  For each $\delta$, we compute two quantities.  The first, shown in Fig.~\ref{fig:boundary-law}(a), is the
normal information $I_{\mathbf n}^K(\delta)$ and
$I_{\mathbf n}^U(\delta)$,
i.e., the averages along the sampled segment of
$\mathbf n^\top C_K(\bx_\delta)\mathbf n$ and
$\mathbf n^\top C_U(\bx_\delta)\mathbf n$, with $C_K$ and $C_U$ defined in
\eqref{eq:CK}--\eqref{eq:CU}.  By Th.~\ref{thm:boundary-loss}, both
should decay like $\delta^2$.  The second, shown in Fig.~\ref{fig:boundary-law}(b), is the spatial median along the
segment of $\sigma_{\min}(\mathsf G(\bx_\delta))$, the smallest singular value
of the increment matrix \eqref{eq:rank3matrix}, which by
Cor.~\ref{cor:rank-boundary} should decay like $\delta$.
Unlike in Sec.~\ref{sec:param_rec}, we compute these quantities from the
gradients of the finite element solutions, and approximate the
Laplacians by finite differences of these gradients.
For the periodic square, we evaluate the same quantities at the same
coordinates.

Fig.~\ref{fig:boundary-law} shows that, on the Neumann square, all three
quantities decrease toward the boundary, whereas the curves for the periodic
square remain approximately flat for small $\delta$.  For $\delta\leq0.1$, the slopes 
are close
to the predicted values of two and one. Moreover, $I_{\mathbf n}^U$ lies below
$I_{\mathbf n}^K$ for every $\delta$.  This holds at every point: since $\bar A$ is the
average of the observed fields, $C_U$ in \eqref{eq:CU} is the covariance of
the gradients, whereas $C_K$ in \eqref{eq:CK} is their second moment about
$\nabla\Ast$, and a covariance never exceeds a second moment about any
other point.  For $\delta>0.1$, the curves are no
longer governed by the boundary: their shape reflects the spatial
structure of the fields, for both the Neumann and the periodic square.

Next, we use a
cubic B-spline parameter space with 24 knot intervals and 48 spatial test intervals in each of $x_1$ and $x_2$, and compute $\lambda_{\min,J}^{1/2}$ as in
Sec.~\ref{sec:param_rec}. To locate where the information about $\eta$ is weakest in the case of unknown $\Ast$, we repeat the
computation with the parameter perturbations restricted to two subspaces:
the span of the basis functions whose supports intersect the boundary strip
$\delta<0.05$ (boundary-related subspace), and the span of the remaining
basis functions (interior subspace).  
Note that, for the boundary-related subspace, the basis function supports may extend beyond this strip.
All three computations use the same
assembled operator, with the static forcing projected out before the
parameter space is restricted.  For the periodic square, there is no
boundary, so the interior subspace coincides with the full space.  We also
compute $u_{\min}\in X_J$, the perturbation of $\eta$ in the full space that attains $\lambda_{\min,J}$, i.e., the one about which the data carry the least information, and the fraction of
its energy $\norm{u_{\min}}_{L^2}^2$ contained in the band $\delta<0.15$. The fraction of energy in this band (which is the same as in Figs.~\ref{fig:square-recovery}
and~\ref{fig:chicago_recovery}) measures how strongly $u_{\min}$ is
localized near the boundary.

\begin{figure}[htb!]
    \centering
    \includegraphics[width=\linewidth]{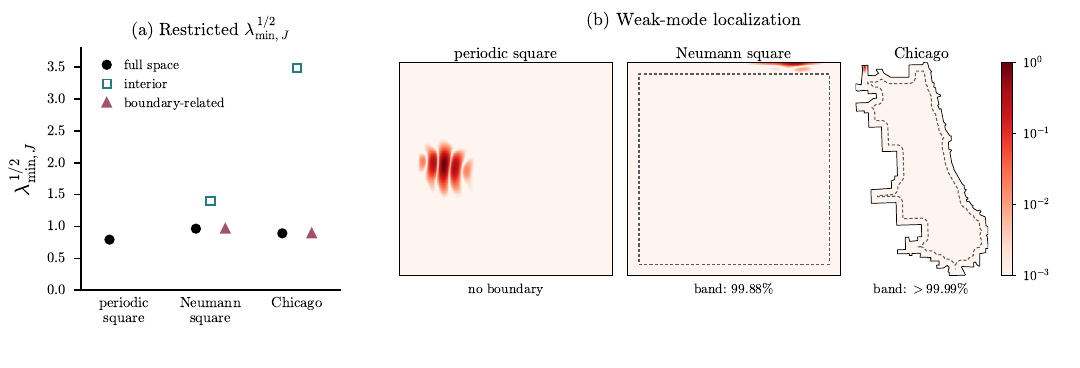}
    \caption{All configurations in the case of unknown $\Ast$:
    (a) $\lambda_{\min,J}^{1/2}$ on the full, interior, and boundary-related
    parameter subspaces and
    (b) energy of $u_{\min}$,
    $|u_{\min}|^2/\max_{\Om}|u_{\min}|^2$, normalized separately in each
    domain.
    Dashed contours mark $\delta=0.15$. The percentages give the fraction of
    the energy of $u_{\min}$ within this band.
}
    \label{fig:weak-modes}
\end{figure}

Fig.~\ref{fig:weak-modes}(a) shows that the full-space values of
$\lambda_{\min,J}^{1/2}$ are comparable in the three configurations. 
We note that for the periodic square the full and interior subspaces coincide and there is no boundary-related value. 
For
the Neumann square and Chicago, restricting the perturbations to the
boundary-related subspace gives nearly the same value as the full space,
whereas restricting them to the interior subspace gives a markedly larger
value.  Hence, the parameter directions that the data constrain least are
associated with the boundary.  Fig.~\ref{fig:weak-modes}(b) confirms this:
$u_{\min}$ has $99.88\%$ of its energy in the boundary band for
the Neumann square and more than $99.99\%$ for Chicago, whereas for the
periodic square it is localized in the interior.  For Chicago, this is
consistent with the concentration near the boundary of the error for
$\eta$ with unknown $\Ast$ observed in Fig.~\ref{fig:chicago-eta-unknown}.

These results agree with the mechanism of Sec.~\ref{sec:boundary-information}: the loss of information
near the boundary is responsible for the smallest eigenvalue of the learner.  

\section{Conclusions}\label{sec:concl}

We presented a weak-form method to learn the spatially varying
parameters of a PDE model of urban crime from observations of the
attractiveness and  criminal density, and we have analyzed what the
data can and cannot determine.  
{The spatially varying parameters are three:
the strength of the neighborhood effects \(\eta\), the static
attractiveness \(\Ast\), and a source term for criminal regeneration \(q\). }
These parameters are identifiable under a
pointwise rank condition on the observed increments of the attractiveness. We gave an exact expression for the information lost when the static attractiveness is unknown.
Near a no-flux boundary, we established the rates at which normal information and the smallest singular value of the rank matrix decay.  Numerical
experiments on a square and on a realistic geometry based on the city of
Chicago showed that the method recovers the parameters with relative errors
of a few percent, and that the predicted losses of information occur.

These results were obtained with noise-free, densely sampled observations
of both the attractiveness and the criminal density, and they prepare the
next steps of this work.  {Part~II will study the effects of observation noise on parameter recovery and extend the model to include police control. } Part~III will consider the case in which only the crime rate is observed, which is closer to the data available in practice.

\section*{Acknowledgments}
Zhong is supported by NSF-AoF grant $\#2225507$.
The authors used Claude (Anthropic) to help edit the text and check the proofs. All content was reviewed and verified by the authors, who take full responsibility for it. The authors declare no competing interests.

\appendix
\section{Prescribed parameters}\label{sec:app}

In this appendix, we give the explicit definitions of the fields
prescribed in the numerical experiments of Sec.~\ref{sec:num_res}.  
To
control their ranges, we use a normalization operator that maps a scalar function $f$ to the interval $[a,b]$:
\[
\mathcal N_{[a,b]}[f]
=
a+(b-a)\frac{f-f_{\min}}{f_{\max}-f_{\min}}.
\]
Here, the
extrema $f_{\min}$ and $f_{\max}$ are evaluated on a fixed reference grid of
spacing $0.01$ over the corresponding domain, independent of the
simulation mesh, so that the prescribed fields do not depend on the
discretization.
As a consequence, the ranges given in Sec.~\ref{sec:num_res} are nominal.  All prescribed analytic fields are
interpolated into the finite element spaces used for the simulations.

The Neumann square and the periodic square share the same prescribed
fields.  Setting $X=2\pi x_1/L$ and $Y=2\pi x_2/L$, with $L=2.8$, we prescribe
\begin{align*}
\eta
&=\mathcal N_{[0.02,0.05]}
 \left[
 \cos(X-3\pi/4)
 +0.8\cos(Y-15\pi/8)
 +0.5\cos(X+Y-\pi)
 \right],\\
q
&=\mathcal N_{[1.15,1.85]}
 \left[
 \cos(X-Y-\pi/2)+0.6\cos(X+Y-\pi/2)
 \right],\\
\Ast
&=\mathcal N_{[1/30,0.19]}
 \left[
 V_{\pi,\pi}+0.6V_{0,\pi}+0.6V_{\pi,0}
 \right],
\end{align*}
where
\[
V_{a,b}(X,Y)
=
\exp\!\left(
0.6[\cos(X-a)-1+\cos(Y-b)-1]
\right).
\]
All three fields are periodic, so the same fields can be used with
periodic boundary conditions on the periodic square.  Moreover, since
$a,b\in\{0,\pi\}$, the derivatives of $V_{a,b}$ normal to each side of the
square vanish, so $\Ast$ satisfies \eqref{eq:nofluxAst} on the Neumann
square.

\begin{figure}[htb!]
    \centering
    \includegraphics[width=0.62\linewidth]
    {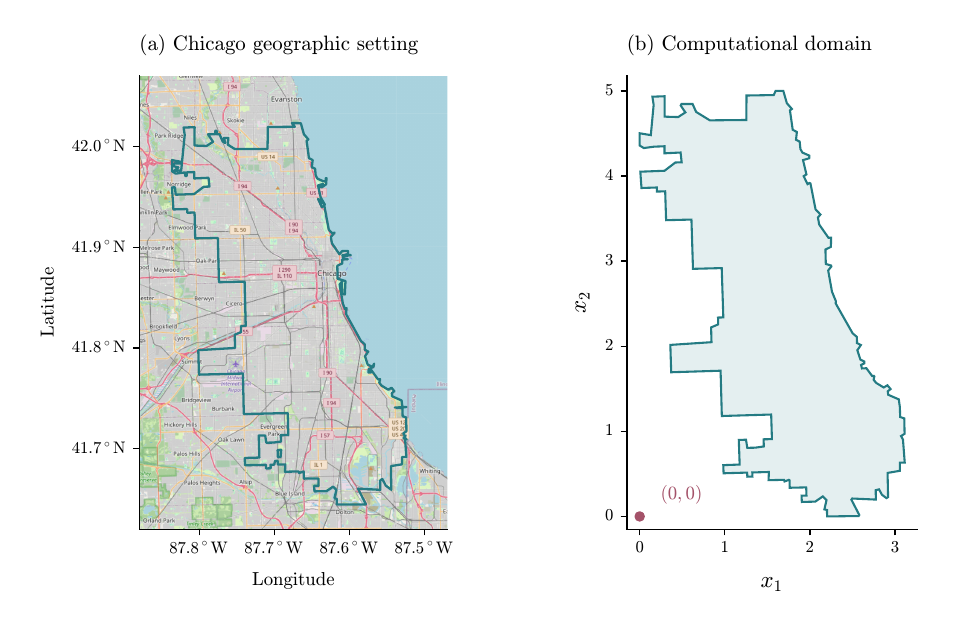}
    \caption{Chicago:
    (a) outline of
    the geographic region used as computational domain and
    (b) computational polygon in the coordinates
    $(x_1,x_2)$ used to define the parameters.
    Basemap from OpenStreetMap
\cite{openstreetmap}.}
    \label{fig:chicago-geometry}
\end{figure}

For the domain of Chicago, Fig.~\ref{fig:chicago-geometry}
shows the geographic setting and the computational coordinates.
The origin is at the lower-left corner of the computational
bounding box, whose longest side is scaled to $5$.
In these coordinates, let
\[
G_{a,b,s}(x_1,x_2)
=
\exp\!\left(-\frac{(x_1-a)^2+(x_2-b)^2}{2s^2}\right)
\]
denote a Gaussian centered at $(a,b)$ with width $s$.
We prescribe
\begin{align*}
\eta
&=\mathcal N_{[0.02,0.05]}
 \left[
 G_{0.95,4.25,0.90}
 -1.2G_{1.65,2.35,0.90}
 +G_{2.45,0.75,0.85}
 \right],\\
q
&=\mathcal N_{[1.15,1.85]}
 \left[
 \frac{x_1-1.55}{1.5\sqrt{2}}
 +\frac{x_2-2.5}{2.5\sqrt{2}}
 +0.8G_{2.40,1.60,0.85}
 \right],\\
S
&=\mathcal N_{[1/30,0.20]}
 \left[
 G_{1.80,3.10,0.75}
 +0.75G_{0.45,4.55,0.75}
 +0.7G_{2.05,0.65,0.75}
 \right].
\end{align*}
The static field $\Ast$ is then obtained from $S$ as described in Sec.~\ref{sec:num_res}.

\begin{thebibliography}{10}

\bibitem{bal2013reconstruction}
G.~Bal and G.~Uhlmann.
\newblock Reconstruction of coefficients in scalar second-order elliptic equations from knowledge of their solutions.
\newblock {\em Commun. Pure Appl. Math.}, 66(10):1629--1652, 2013.

\bibitem{banks1989estimation}
H.~T. Banks and K.~Kunisch.
\newblock {\em Estimation techniques for distributed parameter systems}.
\newblock Birkh{\"a}user, Boston, 1989.

\bibitem{blumstein2002crime}
A.~Blumstein.
\newblock Crime modeling.
\newblock {\em Oper. Res.}, 50(1):16--24, 2002.

\bibitem{mauro2017}
M.~Bongini, M.~Fornasier, M.~Hansen, and M.~Maggioni.
\newblock Inferring interaction rules from observations of evolutive systems {I}: the variational approach.
\newblock {\em Math. Models Methods Appl. Sci.}, 27(5):909--951, 2017.

\bibitem{bortz2023wendy}
D.~M. Bortz, D.~A. Messenger, and V.~Dukic.
\newblock Direct estimation of parameters in {ODE} models using {WENDy}: weak-form estimation of nonlinear dynamics.
\newblock {\em Bull. Math. Biol.}, 85(11):110, 2023.

\bibitem{brantingham1984patterns}
P.~J. Brantingham and P.~L. Brantingham.
\newblock {\em Patterns in crime}.
\newblock Macmillan, New York, 1984.

\bibitem{brunton2016sindy}
S.~L. Brunton, J.~L. Proctor, and J.~N. Kutz.
\newblock Discovering governing equations from data by sparse identification of nonlinear dynamical systems.
\newblock {\em Proc. Natl. Acad. Sci. U.S.A.}, 113(15):3932--3937, 2016.

\bibitem{cai2022multiscale}
Y.~Cai, C.~Wang, and Y.~Zhang.
\newblock A multiscale stochastic criminal behavior model and the convergence to a piecewise-deterministic-{Markov}-process limit.
\newblock {\em Math. Models Methods Appl. Sci.}, 32(4):619--645, 2022.

\bibitem{calamunci2022highway}
F.~Calamunci and J.~Lonsky.
\newblock The road to crime: an unintended consequence of the {Interstate Highway System}.
\newblock {\em Econ. J.}, 135(667):748--772, 2025.

\bibitem{chen2026kellersegel}
C.~Chen, C.~Liu, and M.~Zhong.
\newblock Learning profile functions in {Keller--Segel} models.
\newblock {\em Adv. Comput. Sci. Eng.}, 10:21--46, 2026.

\bibitem{Chicago}
{City of Chicago}.
\newblock Boundaries -- city -- map.
\newblock Chicago Data Portal, \url{https://data.cityofchicago.org/d/ewy2-6yfk}.

\bibitem{dorsogna2015statistical}
M.~R. D'Orsogna and M.~Perc.
\newblock Statistical physics of crime: a review.
\newblock {\em Phys. Life Rev.}, 12:1--21, 2015.

\bibitem{ersoz2025artificial}
F.~Ers{\"o}z, T.~Ers{\"o}z, F.~Marcelloni, and F.~Ruffini.
\newblock Artificial intelligence in crime prediction: a survey with a focus on explainability.
\newblock {\em IEEE Access}, 13:59646--59674, 2025.

\bibitem{falk1983error}
R.~S. Falk.
\newblock Error estimates for the numerical identification of a variable coefficient.
\newblock {\em Math. Comput.}, 40(162):537--546, 1983.

\bibitem{feng2024learning}
J.~Feng and M.~Zhong.
\newblock Learning collective behaviors from observation.
\newblock In S.~Foucart and S.~Wojtowytsch, editors, {\em Explorations in the Mathematics of Data Science: The Inaugural Volume of the Center for Approximation and Mathematical Data Analytics}, pages 101--132. Birkh{\"a}user, 2024.

\bibitem{gao2026self}
Y.~Gao, Q.~Lang, and F.~Lu.
\newblock Self-test loss functions for learning weak-form operators and gradient flows.
\newblock {\em J. Sci. Comput.}, 108(1):28, 2026.

\bibitem{geuzaine2009gmsh}
C.~Geuzaine and J.-F. Remacle.
\newblock {Gmsh}: a 3-d finite element mesh generator with built-in pre- and post-processing facilities.
\newblock {\em Int. J. Numer. Methods Eng.}, 79(11):1309--1331, 2009.

\bibitem{guo2024sde}
Z.~Guo, I.~Cialenco, and M.~Zhong.
\newblock Learning stochastic dynamics from data.
\newblock In {\em ICLR Workshop on AI4DifferentialEquations in Science}, 2024.

\bibitem{gurevich2019robust}
D.~R. Gurevich, P.~A.~K. Reinbold, and R.~O. Grigoriev.
\newblock Robust and optimal sparse regression for nonlinear {PDE} models.
\newblock {\em Chaos}, 29(10):103113, 2019.

\bibitem{OurCode}
B.~Hao, K.~Mily, A.~Quaini, and M.~Zhong.
\newblock {UrbanCrime-Sim}.
\newblock Computer software, \url{https://github.com/baolihao/UrbanCrime-Sim}.

\bibitem{Hao2026b}
B.~Hao, K.~Mily, A.~Quaini, and M.~Zhong.
\newblock Crime hotspot dynamics in residential burglary models with police response.
\newblock arXiv:2605.17709, 2026.

\bibitem{Hao2026}
B.~Hao, K.~Mily, A.~Quaini, and M.~Zhong.
\newblock A finite element framework for simulating residential burglary in realistic urban geometries.
\newblock {\em Math. Models Methods Appl. Sci.}, 36(5):1019--1050, 2026.

\bibitem{johnson1997new}
S.~D. Johnson, K.~Bowers, and A.~Hirschfield.
\newblock New insights into the spatial and temporal distribution of repeat victimization.
\newblock {\em Br. J. Criminol.}, 37(2):224--241, 1997.

\bibitem{jones2010statistical}
P.~A. Jones, P.~J. Brantingham, and L.~R. Chayes.
\newblock Statistical models of criminal behavior: the effects of law enforcement actions.
\newblock {\em Math. Models Methods Appl. Sci.}, 20(Suppl. 1):1397--1423, 2010.

\bibitem{kang2021ident}
S.~H. Kang, W.~Liao, and Y.~Liu.
\newblock {IDENT}: identifying differential equations with numerical time evolution.
\newblock {\em J. Sci. Comput.}, 87(1):1, 2021.

\bibitem{kohn1988variational}
R.~V. Kohn and B.~D. Lowe.
\newblock A variational method for parameter identification.
\newblock {\em ESAIM Math. Model. Numer. Anal.}, 22(1):119--158, 1988.

\bibitem{lancaster1981surfaces}
P.~Lancaster and K.~Salkauskas.
\newblock Surfaces generated by moving least squares methods.
\newblock {\em Math. Comput.}, 37(155):141--158, 1981.

\bibitem{lang2022meanfield}
Q.~Lang and F.~Lu.
\newblock Learning interaction kernels in mean-field equations of first-order systems of interacting particles.
\newblock {\em SIAM J. Sci. Comput.}, 44(1):A260--A285, 2022.

\bibitem{lloyd2016exploring}
D.~J.~B. Lloyd, N.~Santitissadeekorn, and M.~B. Short.
\newblock Exploring data assimilation and forecasting issues for an urban crime model.
\newblock {\em Eur. J. Appl. Math.}, 27(3):451--478, 2016.

\bibitem{lu2019pnas}
F.~Lu, M.~Zhong, S.~Tang, and M.~Maggioni.
\newblock Nonparametric inference of interaction laws in systems of agents from trajectory data.
\newblock {\em Proc. Natl. Acad. Sci. U.S.A.}, 116(29):14424--14433, 2019.

\bibitem{mei2020existence}
L.~Mei and J.~Wei.
\newblock The existence and stability of spike solutions for a chemotaxis system modeling crime pattern formation.
\newblock {\em Math. Models Methods Appl. Sci.}, 30(9):1727--1764, 2020.

\bibitem{messenger2021wsindypde}
D.~A. Messenger and D.~M. Bortz.
\newblock Weak {SINDy} for partial differential equations.
\newblock {\em J. Comput. Phys.}, 443:110525, 2021.

\bibitem{messenger2021wsindymms}
D.~A. Messenger and D.~M. Bortz.
\newblock Weak {SINDy}: {Galerkin}-based data-driven model selection.
\newblock {\em Multiscale Model. Simul.}, 19(3):1474--1497, 2021.

\bibitem{messenger2022meanfield}
D.~A. Messenger and D.~M. Bortz.
\newblock Learning mean-field equations from particle data using {WSINDy}.
\newblock {\em Physica D}, 439:133406, 2022.

\bibitem{mohler2020impact}
G.~Mohler, A.~L. Bertozzi, J.~Carter, M.~B. Short, D.~Sledge, G.~E. Tita, C.~D. Uchida, and P.~J. Brantingham.
\newblock Impact of social distancing during {COVID}-19 pandemic on crime in {Los} {Angeles} and {Indianapolis}.
\newblock {\em J. Crim. Justice}, 68:101692, 2020.

\bibitem{nurhan2026modelling}
Y.~Nurhan and M.~B. Short.
\newblock Modelling crime response to deterrence: existence of solutions, optimal policies, and fairness.
\newblock {\em Eur. J. Appl. Math.}, 37(5):1222--1243, 2026.

\bibitem{openstreetmap}
{OpenStreetMap contributors}.
\newblock {OpenStreetMap}.
\newblock \url{https://www.openstreetmap.org/copyright}.

\bibitem{pan2018crime}
C.~Pan, B.~Li, C.~Wang, Y.~Zhang, N.~Geldner, L.~Wang, and A.~L. Bertozzi.
\newblock Crime modeling with truncated {L{\'e}vy} flights for residential burglary models.
\newblock {\em Math. Models Methods Appl. Sci.}, 28(9):1857--1880, 2018.

\bibitem{park2021investigating}
J.~Park, F.~P. Schoenberg, A.~L. Bertozzi, and P.~J. Brantingham.
\newblock Investigating clustering and violence interruption in gang-related violent crime data using spatial--temporal point processes with covariates.
\newblock {\em J. Am. Stat. Assoc.}, 116(536):1674--1687, 2021.

\bibitem{reinbold2020noisy}
P.~A.~K. Reinbold, D.~R. Gurevich, and R.~O. Grigoriev.
\newblock Using noisy or incomplete data to discover models of spatiotemporal dynamics.
\newblock {\em Phys. Rev. E}, 101(1):010203(R), 2020.

\bibitem{richter1981inverse}
G.~R. Richter.
\newblock An inverse problem for the steady state diffusion equation.
\newblock {\em SIAM J. Appl. Math.}, 41(2):210--221, 1981.

\bibitem{rodriguez2010local}
N.~Rodriguez and A.~L. Bertozzi.
\newblock Local existence and uniqueness of solutions to a {PDE} model for criminal behavior.
\newblock {\em Math. Models Methods Appl. Sci.}, 20(Suppl. 1):1425--1457, 2010.

\bibitem{rodriguez2020relaxation}
N.~Rodr{\'i}guez and M.~Winkler.
\newblock Relaxation by nonlinear diffusion enhancement in a two-dimensional cross-diffusion model for urban crime propagation.
\newblock {\em Math. Models Methods Appl. Sci.}, 30(11):2105--2137, 2020.

\bibitem{rudy2019parametric}
S.~Rudy, A.~Alla, S.~L. Brunton, and J.~N. Kutz.
\newblock Data-driven identification of parametric partial differential equations.
\newblock {\em SIAM J. Appl. Dyn. Syst.}, 18(2):643--660, 2019.

\bibitem{rudy2017pdefind}
S.~H. Rudy, S.~L. Brunton, J.~L. Proctor, and J.~N. Kutz.
\newblock Data-driven discovery of partial differential equations.
\newblock {\em Sci. Adv.}, 3(4):e1602614, 2017.

\bibitem{santitissadeekorn2018sequential}
N.~Santitissadeekorn, M.~B. Short, and D.~J.~B. Lloyd.
\newblock Sequential data assimilation for {1D} self-exciting processes with application to urban crime data.
\newblock {\em Comput. Stat. Data Anal.}, 128:163--183, 2018.

\bibitem{short2010nonlinear}
M.~B. Short, A.~L. Bertozzi, and P.~J. Brantingham.
\newblock Nonlinear patterns in urban crime: hotspots, bifurcations, and suppression.
\newblock {\em SIAM J. Appl. Dyn. Syst.}, 9(2):462--483, 2010.

\bibitem{short2010dissipation}
M.~B. Short, P.~J. Brantingham, A.~L. Bertozzi, and G.~E. Tita.
\newblock Dissipation and displacement of hotspots in reaction-diffusion models of crime.
\newblock {\em Proc. Natl. Acad. Sci. U.S.A.}, 107(9):3961--3965, 2010.

\bibitem{short2008statistical}
M.~B. Short, M.~R. D'Orsogna, V.~B. Pasour, G.~E. Tita, P.~J. Brantingham, A.~L. Bertozzi, and L.~B. Chayes.
\newblock A statistical model of criminal behavior.
\newblock {\em Math. Models Methods Appl. Sci.}, 18(Suppl. 1):1249--1267, 2008.

\bibitem{stuart2010bayesian}
A.~M. Stuart.
\newblock Inverse problems: a {B}ayesian perspective.
\newblock {\em Acta Numer.}, 19:451--559, 2010.

\bibitem{tang2025wgident}
C.~Tang, R.~Y. He, and H.~Liu.
\newblock {WG-IDENT}: weak group identification of {PDEs} with varying coefficients.
\newblock {\em J. Comput. Phys.}, 545:114454, 2026.

\bibitem{tang2023weakident}
M.~Tang, W.~Liao, R.~Kuske, and S.~H. Kang.
\newblock {WeakIdent}: weak formulation for identifying differential equation using narrow-fit and trimming.
\newblock {\em J. Comput. Phys.}, 483:112069, 2023.

\bibitem{tran2026weak}
A.~Tran and D.~M. Bortz.
\newblock Weak form scientific machine learning: test function construction for system identification.
\newblock {\em SIAM J. Sci. Comput.}, 48(4):C890--C915, 2026.

\bibitem{wang2019deep}
B.~Wang, P.~Yin, A.~L. Bertozzi, P.~J. Brantingham, S.~J. Osher, and J.~Xin.
\newblock Deep learning for real-time crime forecasting and its ternarization.
\newblock {\em Chin. Ann. Math. Ser. B}, 40(6):949--966, 2019.

\bibitem{wang2017deep}
B.~Wang, D.~Zhang, D.~Zhang, P.~J. Brantingham, and A.~L. Bertozzi.
\newblock Deep learning for real time crime forecasting.
\newblock In {\em Proceedings of the 2017 International Symposium on Nonlinear Theory and Its Applications (NOLTA)}, pages 330--333, 2017.

\bibitem{weisburd2015law}
D.~Weisburd.
\newblock The law of crime concentration and the criminology of place.
\newblock {\em Criminology}, 53(2):133--157, 2015.

\bibitem{wilson1982broken}
J.~Q. Wilson and G.~L. Kelling.
\newblock Broken windows: the police and neighborhood safety.
\newblock {\em Atl. Mon.}, 249(3):29--38, 1982.

\bibitem{woodworth2014nonlocal}
J.~T. Woodworth, G.~O. Mohler, A.~L. Bertozzi, and P.~J. Brantingham.
\newblock Non-local crime density estimation incorporating housing information.
\newblock {\em Philos. Trans. R. Soc. A}, 372:20130403, 2014.

\bibitem{wright1994burglars}
R.~T. Wright and S.~H. Decker.
\newblock {\em Burglars on the job: streetlife and residential break-ins}.
\newblock Northeastern Univ. Press, Boston, 1994.

\bibitem{zipkin2014cops}
J.~R. Zipkin, M.~B. Short, and A.~L. Bertozzi.
\newblock Cops on the dots in a mathematical model of urban crime and police response.
\newblock {\em Discrete Contin. Dyn. Syst. Ser. B}, 19(5):1479--1506, 2014.

\end{thebibliography}
\end{document}